\documentclass{statsoc}
\pdfoutput=1 

\usepackage{amsmath}
\usepackage{amssymb}
\usepackage{amsfonts}
\usepackage{dsfont}
\usepackage{epstopdf}
\usepackage{graphicx}
\usepackage{fixltx2e} 
\usepackage{xcolor}
\usepackage{hyperref}
\hypersetup{
    colorlinks=true,
    linkcolor=blue, 
    citecolor=blue,
    urlcolor=blue,
    pdftitle={Estimating the Parameters of the Waxman Graph},
    breaklinks
}

\definecolor{palepink}{rgb}{1.0,0.9,0.9}
\usepackage{natbib}
\usepackage[a4paper]{geometry}

\usepackage{framed}
\usepackage{lipsum}
\colorlet{shadecolor}{blue!20}

\usepackage[utf8]{inputenc}
\usepackage{newunicodechar}
\newunicodechar{✓}{\checkmark}

\usepackage{fixltx2e} 

\usepackage{epstopdf}

\newlength{\figurewidthA}
\newlength{\figurewidthB}
\newlength{\captionspace}
\newlength{\figstarspace}
\newcommand{\ie}{{\em i.e., }}
\newcommand{\eg}{{\em e.g., }}
\newcommand{\etal}{{\em et al.}}

\newcommand{\anonymize}[1]{anonymous}

\newtheorem{thm}{Theorem}

\newcounter{subfig_counter}
\newenvironment{subfigure}[2][t]%
{\begin{minipage}[#1]{#2}%
    \renewcommand{\caption}[1]{##1}%
}%
{\end{minipage}}%

\def\equationautorefname~#1\null{%
  (#1)\null
}

\newcommand{\Prob}[1]{\mathds{P} \! \left\{#1\right\}}
\newcommand{\E}[1]{\mathds{E} \! \left[#1\right]}

\title{Estimating the Parameters of the Waxman Random Graph}

\author{Matthew Roughan}
\address{ARC Centre of Excellence for Mathematical \& Statistical
         Frontiers,
         School of Mathematical Sciences, University of Adelaide, Australia.}
\email{matthew.roughan@adelaide.edu.au}

\author{Jonathan Tuke}
\address{School of Mathematical Sciences, University of Adelaide, Australia.}
\email{simon.tuke@adelaide.edu.au}

\author[Roughan {\em et al.}]{Eric Parsonage}
\address{School of Mathematical Sciences, University of Adelaide, Australia.}
\email{eric.parsonage@adelaide.edu.au}

\begin{document}

\maketitle

\begin{abstract}
  The Waxman random graph is a generalisation of the simple
  Erd\H{o}s-R\'enyi or Gilbert random graph. It is useful for modelling  physical networks 
  where the increased cost of longer links means they are less likely to be built, and thus less numerous than shorter links. 
  The model has been in continuous use for over two decades with many attempts to select parameters which match real
  networks. In most the parameters have been arbitrarily selected, but  there are a few cases where they have been
  calculated using a formal estimator. However, the performance of the estimator was not evaluated in any of these cases. 
  This paper presents both the first evaluation of formal estimators for the parameters of these graphs, and a new Maximum
  Likelihood Estimator with $O(n)$ computational time complexity that requires only link lengths as input.
\end{abstract}

\section{Introduction}
\label{sec:intro}

The study of random graphs provides insight into
the formation of real networks and methods to synthesise test
networks for use in simulations. 
Early research was dominated by the study of
Gilbert-Erd\H{o}s-R\'enyi (GER) random graphs
\citep{gilbert59:_random_graph,erdos60}, which link two vertices
independently with a fixed probability, and hence select particular
graphs (with a fixed numbers of nodes) all with equal probability.
Due partly to its mathematical tractability this model was studied for 
many years despite clear limitations in its applicability to real networks.

Random graphs explicitly incorporating geometry became popular with
the introduction of the Waxman random graph
\citep{b.m.waxman88:_graph}, proposed as an alternative to the GER
random graph as a more realistic setting for testing networking
algorithms.  The Waxman graph has subsequently been used in
applications as wide-ranging as computer networks, transportation and
biology (Waxman's original paper is listed by Google Scholar as having
been cited around 3000 times). Moreover, the model has been reinvented
at least once so not all use cases can be traced back to this paper.

The GER random graph links every pair of vertices independently with a
fixed probability, whereas the Waxman graph reflects that in real
networks longer links are often more costly or difficult to construct,
and their existence therefore less likely. It links nodes $i$ and $j$
with a probability given by a function of the distance $d_{ij}$
between them. The form chosen by Waxman was the negative
exponential\footnote{Note that our parameterisation differs from that
  in much of the literature, for reasons to be explained in the
  following sections.}, \ie
\begin{equation}
  \label{eq:waxman_prob}
  p(d_{ij}) = q \, e^{- s d_{ij}},
\end{equation}
for parameters $q,s \geq 0$.  A Waxman random graph is generated by
randomly choosing a set of points in a section of the plane (usually
the unit square), and then linking these points independently
according to the distance between them and (\ref{eq:waxman_prob}).

Despite frequent use of the Waxman graph there is no formal
literature on how to estimate the parameters $(q,s)$ from a given
graph, or set of graphs. In most works the parameter values have been
chosen arbitrarily with the authors giving little or no justification for the
values used. In many cases authors use the parameter
values given in earlier works without regard for their origin or
applicability. A few works use more careful estimates, 
but do not evaluate the estimator performance. They
simply apply estimates and report results.

This paper presents the first work of which we are aware on the formal
estimation of the parameters of the Waxman graph.  As part of our
investigation we explain many properties of the
Waxman graph and its likelihood function, and derive:
\begin{itemize}
\item its average node degree, expected number of edges, and
  average edge length;
\item asymptotic expressions for these for large $s$; 
\item  the Kullback-Leibler divergence between Waxman graphs
  and the GER graph;
\item  minimal sufficient statistics for the model and
  hence a MLE for the model parameters. 
\end{itemize}

We present the Maximum Likelihood Estimator (MLE) under the assumption
of independence between links, and demonstrate that its performance is
close to the Cram\'{e}r-Rao lower bound. We also compare the MLE to
other existing approaches and show its advantages in:
\begin{itemize}
\item accuracy;
\item computational time complexity (by suitable sampling of the existing link lengths $O(n)$ can always be achieved);
\item minimal input (using only the lengths of links which exist or a sufficiently large sample of such links). 
\end{itemize}

Finally we use the MLE  to estimate the distance dependence of three real networks, using one biological and two Internet datasets. 


\section{Background and Related Work}
\label{sec:background}

A graph (or network) consists of a set of $n$ vertices (synonymously
referred to as nodes) which without loss of generality we label ${\cal
  V}=\{1,2,\ldots, n\}$, and a set of edges (or links) ${\cal E}
\subset {\cal V} \times {\cal V}$. We denote the number of edges by $e
= |{\cal E}|$. Here we are primarily concerned with undirected graphs,
though much work on random graphs is easily generalised to directed
graphs.

Two nodes $i$ and $j$ are {\em adjacent} or {\em
  neighbours} if $(i,j) \in {\cal E}$. They are {\em
  connected} if a path $i=n_1\!-\!n_2\!-\cdots\!-\!n_k=j$ exists, where
$(n_l,n_{l+1}) \in {\cal E}$, and the graph is connected if
all pairs of nodes $(i,j)$ are connected.
 
The GER random graph $G_{n,p}$ of $n$ vertices is constructed by
assigning independently and with a fixed probability each tuple of
nodes $(i,j)$ to be in ${\cal E}$.  The Waxman random graph
generalises this by making the probability of each tuple of nodes an
edge dependent on a function $p$ of the distance between the two
nodes\footnote{Some examples of use of the Waxman graph
  include~\citep{b.m.waxman88:_graph,thomas94:_waxman,zegura96:_gtitm,zegura97:_quant_compar_of_graph_based,doar_leslie,wei94:_trade_offs_multic_trees_algor,Salama97evaluationof,verma98:_qos,matta99:_qdmr,shaikh99:_load_sensit_routin_long_lived_ip_flows,fortz00:or_ospf,neve00:_tamcr,Rastogi01algorithmsfor,Wu200029,guo03:_searc_qos,kuipers04:_qualit,kaiser04:_spatial,Kaiser2004297,gunduz04:_cell_graph_of_cancer,Carzaniga04arouting,Lua05onthe,holzer05:_combin_speed_up_techn_for,Wang05cross-layeroptimization,ahuja09:_singl_link_failur_detec_in,huang07:_multic_routin_based_ant_system,malladi07:_improv_secur_commun_polic_agreem_build_coalit,tran09:_pub_sub,fang11:_iterat,fang11,costa10:_model_evolut_compl_networ_throug,janssen10:_spatial_model_virtual_networ,davis14:_spatial}.}.

While the Waxman graph was originally defined using Euclidean 
distances on a rectangle or straight line with points on an integer grid, and most
subsequent work has considered graphs defined with points randomly placed in
the unit square, there is no impediment to choosing points in
an arbitrary convex region\footnote{Convexity is not strictly
  required, except where there is the notion that the links are
  physical, and must themselves lie in the region of interest.} with an
arbitrary distance metric. Hence, we define a Waxman graph by placing
$n$ nodes uniformly at random within some defined region $R$ of a
metric space $\Omega$ with a distance metric $d(x,y)$ and each pair of nodes is made adjacent independently with probability
given by (\ref{eq:waxman_prob}).

 This differs from Waxman's
original parameterisation \footnote{Although users of the Waxman model commonly state
  that $\alpha \in (0,1]$ there is no reason it should not take
  arbitrary positive values.}~\citep{b.m.waxman88:_graph}
\begin{equation}
  \label{eq:waxman_orig}
  p(d_{ij}) = \beta \, e^{- d_{ij}/L \alpha}
\end{equation}
for parameters $\alpha, \beta > 0$, where $L$ is the greatest distance
possible between any two points in the region of interest. The advantage of
Waxman's parameterisation is that $\alpha$ is a dimensionless
constant. Unfortunately, previous authors have confused this notation by reversing
the roles of $\alpha$ and $\beta$ with almost equal regularity, to the
point where  parameters chosen in one paper have been reversed in another 
purporting to compare results. Hence, we provide an alternate parameterisation 
\autoref{eq:waxman_prob} chosen with the estimation problem in mind. 
We deliberately avoid a dimensionless parameter because real estimates naturally have
units and it has been common in the literature to assume incorrectly that use of such a parameter removes dependence on the region shape.

The parameter $s$ determines how dependent links are on
distance. Larger $s$ decreases the likelihood of long links.

The use of $\beta$ in past works on the Waxman graphs has allowed for
$\beta > 1$, at which point the function \autoref{eq:waxman_orig} must
be truncated, creating a corner in the deterrence function
shape. This is undesirable from the point of view of parametric
estimation, so we restrict $q \in (0,1]$. In this case $q$ is a
thinning of the edges, and we shall see that this allows decoupling of
estimation of the two parameters.

There are many variants of the Waxman graph through modifications of
the distance deterrence function, for instance
\citet{zegura97:_quant_compar_of_graph_based} proposed $p(d; \beta) =
\beta \exp(-d/(L - d))$. One of the more interesting variants, which
we refer to here as the DASTB (Davis-Abbasi-Shah-Telfer-Begon) graph
\citep{davis14:_spatial}, takes Waxman-like rates for contacts between
the agents represented by vertices in the graph, and then considers an
edge to exist if at least one contact is made. The number of contacts
is assumed to have Poisson distribution with parameter $k_{ij}$ given
by the Waxman function \autoref{eq:waxman_prob}, and then the
probability of at least one contact is $p_{ij}(d) = 1 - \exp( -k_{ij}
)$. We do not formally consider this model here (though we do adapt
their estimation approach) except to note that for sparse networks,
the low probability of contact means that the two ensembles of random
graphs are similar.

Properties such as connectivity and path lengths in Waxman graphs have
been studied in several
works~\citep{thomas94:_waxman,zegura96:_gtitm,zegura97:_quant_compar_of_graph_based,Mieghem01,m.naldi05:_connec_of_waxman_graph}. However, few have considered the estimation of parameters. Indeed, some of the earlier works compared graphs merely by visual inspection.

The first work we are aware of which attempted to estimate the
parameters of a Waxman graph from real data is that of
\citet{Lakhina:2002:GLI:637201.637240}.  Using a large set of real
Internet and geographic data the authors found that an exponential
distance-based probability was reasonable.  The authors conducted the study
with some care, comparing two sets of data and finding consistent results by using a log-linear regression of the
link distance function against the link distances.
However they made no effort to consider the efficiency or accuracy of their estimator.
 
 
Estimation has also been used on the DASTB
model~\citep{davis14:_spatial} by determining the parameters of a
binary Generalised Linear Model (GLM). We can adapt this approach to
the almost equivalent Waxman graphs by considering that links are the
binary outcomes of a treatment which is functionally  dependent on link
distances through the distance deterrence function, then use a GLM to estimate the parameters.

\citet{davis14:_spatial} showed, using their parameter estimations and
Akaike's Information Criterion corrected for small sample sizes (AICc),
that the DASTB graph was a better model for their data than the GER
(and provided further comparisons against more complex models), and
showed a correlation between $s$ and the population density. 
This has implications for pathogen persistence in
wildlife because density has often been considered a critical
parameter of disease transmission. However, if contact over longer
distances increases in lower-density populations then pathogen
persistence may not be so sensitive to host density, as dense
populations do not mix so freely.

The underlying assumption of \citet{Lakhina:2002:GLI:637201.637240}
and \citet{davis14:_spatial} is that {\em all} of the distances are
known, even for links that do not exist. This is possible if all node
locations are known, but in practice it is sometimes easier to measure
the length of a link that exists -- say by probing it -- than to
estimate the length of a hypothetical link. Consider also graphs for
which the ``distance'' is not a physical quantity, but rather a cost
or an administratively configured link-weight, for instance in
computer network routing protocols such as OSPF link weights can be
the output of an optimisation \citep{fortz03:inoc}. In these cases
link ``distance'' does not even exist for non-links. 

We shall also show that techniques which depend on existence tests of
all edges have time complexity $O(n^2)$, where $n$ is the number of
nodes in the graph, whereas our estimation method has time complexity
$O(e)$, where $e$ is the number of edges in the graph. This can
results in a dramatic improvement in computation time as large real
graphs are often very sparse.

Waxman graphs are an example of the general class of SERNs (Spatially
Embedded Random Graphs) \citep{kosmidis08:_struc}. One aim of this
work is to develop intuition which can be extended to estimation of
the parameters of random graphs from the general class. We will leave
consideration of the general case to later work because complex
questions of existence and identifiability arise.

\section{General properties of Waxman graphs}

The starting point for the creation of a Waxman random graph is to
generate a set of $n$ points uniformly at random in some region $R$ of
a metric space $\Omega$. For any given region we can derive a
probability density function $g(t)$ for the distance between an
arbitrary pair of random points. This is commonly called the
Line-Picking-Problem: for instances of regions including lines, balls,
spheres, cubes, and rectangles see
\citet{b.ghosh51:_random_rect,Rosenberg200499}
and
\citet{weisstein:_line_line_picking,weisstein:_circle_line_picking,weisstein:_squar_line_picking,weisstein:_sphere_line_picking,weisstein:_ball_line_picking}.

Given the distribution of distances between points, we can calculate
the probability that an arbitrary link exists (prior to knowing the
distances):
\begin{equation}
  \label{eq:link_prob1}
   \Prob{ (i,j) \in {\cal E} \mid q, s } = q \int_0^\infty \exp(- s t) g(t)\, dt = q \tilde{G}(s) ,
\end{equation}
for any $i \neq j$, where $\tilde{G}(s)$ is the Laplace transform of
$g(t)$ (or equivalently it is the moment generating function w.r.t. to
$-s$). We know that the Laplace transform at $s=0$ of a probability
density is the normalisation constraint, so $\tilde{G}(0)=1$. Hence
when $s=0$ there is no distance dependence and the result is the GER random
graph.

From this probability we can also compute features of the graph such
as the average node degree
\begin{equation}
  \label{eq:node_degree}
  \bar{k} = (n-1) q \tilde{G}(s),
\end{equation}
from which we can derive values of $q$ that produce given average
degree for a given network size and $s$. From the handshake theorem we
can derive the average number of edges to be
\begin{equation}
  \label{eq:average_edges}
   \bar{e} = n (n-1) q \tilde{G}(s)/ 2.
\end{equation}

We can then derive the distribution of the length $d$ of a link in the
Waxman graph, and we denote this by $f( d \, | q, s) = \Prob{d_{ij} =
  d \mid {(i,j)} \in {\cal E} }$ which is
\begin{eqnarray}
f(d \mid q, s)
  & = & \frac{\Prob{{(i,j)} \in {\cal E} \mid d_{(i,j)} = d; q, s } \Prob{d_{ij} = d } }{\Prob{{(i,j)} \in {\cal E} \mid q, s}} \nonumber \\
  & = & \frac{g(d) \exp(- s d)}{\tilde{G}(s)} .
\label{eqn:link_len_dist}
\end{eqnarray}
Note the $q$ is a thinning parameter, and thus should not change the
link length distribution, and we can see that it vanishes from the
distribution. Hence we generally write $f$ omitting $q$, \ie $f(d \mid
s)$.

The formula allows easy calculation of the mean length of links in the
Waxman graph
\begin{equation}
 \E{d \mid s}  =  \frac{1}{\tilde{G}(s)} \int_{0}^{\infty} t g(t) \exp(- s t) \, dt 
  =  -\frac{\tilde{G}'(s)}{\tilde{G}(s)}.
\label{eqn:av_link_dist}
\end{equation}

For small distances $t$, region-boundary effects are minimal, and so
the function $g(t)$ depends only on the dimension of the embedding
space. For example the square and disk both have $g(t) \stackrel{t
  \downarrow 0}{\sim} 2 \pi t$, which comes from the size of the ring
of radius $t$. Given a Euclidean distance metric on a $k$-dimensional
space the small $t$ approximation is
\begin{equation}
  g(t) \simeq \frac{k  \pi^{k/2}}{\Gamma(k/2+1)} t^{k-1},
\end{equation}
which is the surface area of the $(k-1)$-sphere (the surface of the
hyper-sphere embedded in the $k$-dimensional Euclidean space).

Using the following theorem, we can derive large $s$ approximations
for the Laplace transforms of a Cumulative Density Function (CDF). We
work here mainly with the probability density function $g(t)$, and
this must be integrated to get the CDF, but otherwise the
applicability should be clear.

\begin{thm}[\citet{Feller_2}, pp.445-6]
 \label{theorem1}
  If $0 < a < \infty$, and we have a positive measure concentrated on
  $(0,\infty)$ defined by the CDF $H(t)$ with Laplace-Stieltjes
  transform $\tilde{H}(s)$, then
  \begin{equation}
    \label{eq:tauberian}
      H(t) \stackrel{t \rightarrow 0}{\sim} \frac{L(t)}{\Gamma(a+1)} t^a
     \;\; \Leftrightarrow \;\;
      \tilde{H}(s) \stackrel{s \rightarrow \infty}{\sim} L(1/t) t^{-a},
  \end{equation}
  where $L(t)$ is a {\em slowly varying} function at $0$, where this
  is defined as a function where $\lim_{t \rightarrow 0}
  L(xt)/L(t) = 1$ for all $x >0$.
\end{thm}

Applying \autoref{theorem1} for a $k$-dimensional Euclidean space we
get
\begin{eqnarray}
   \tilde{G}(s)  & \stackrel{s \rightarrow \infty}{\sim} &  \frac{\pi^{k/2} \Gamma(k+1)  }{\Gamma(k/2+1)} s^{-k},  \label{eq:large_s_Gs}\\
   \tilde{G}'(s) & \stackrel{s \rightarrow \infty}{\sim} &  \frac{-n\pi^{k/2} \Gamma(k+1)}{\Gamma(k/2+1)} s^{-k-1}, \label{eq:large_s_Gsd}\\
   \E{d \mid s}  & \stackrel{s \rightarrow \infty}{\sim} & k / s.  \label{eq:large_s_d}
\end{eqnarray}

The Waxman graph is usually defined in terms of a number of nodes and
two parameters $(q,s)$. However, the critical feature relating these
two is the node density on the space $\rho$.  This parameter
determines the number of nodes within a certain distance of a given
node.  \citet{davis14:_spatial} showed the crucial importance of the
relationship between density and the spatial deterrence function. We
can see from the small $t$, large $s$ approximations that this
relationship is (absent of boundary effects) a simple inverse
relationship that generalises to arbitrary dimension. Intuitively it
arises because for large $s$ only short links are likely, and hence
the boundaries of the region play little part in the distribution.

As a measure of the difference between the Waxman and GER link-distance distribution we use the
Kullback-Leibler Divergence (KLD)\footnote{The KLD is not strictly a
  distance metric because it lacks symmetry, hence the term
  ``divergence''.}, which provides an Information Theoretic measure of
the difference.

The KLD of distribution $f(t|s)$ from $g(t)$ is defined by
\begin{equation}
  KLD\big( g(t) \big| \! \big| f(t|s) \big)
  = \int_0^\infty g(t) \ln \frac{g(t)}{f(t|s)} \, dt .
\end{equation}
From \autoref{eqn:link_len_dist} we can derive
\begin{eqnarray}
  KLD\big( g(t) \big| \! \big| f(t|s) \big)
  & = & \int_0^\infty g(t) \big[ \ln g(t) - \ln f(t|s) \big] \, dt \nonumber \\
  & = & \int_0^\infty g(t) \big[ \ln g(t) - \ln g(t) + st + \ln \tilde{G}(s) \big] \, dt \nonumber \\
  & = & s \int_0^\infty t g(t)\, dt + \ln \tilde{G}(s) \int_0^\infty g(t)\, dt \nonumber \\
  & = & s \E{d} + \ln \tilde{G}(s).
  \label{eq:kld}
\end{eqnarray}
Taking the small $t$, large $s$ approximations we can see that
\begin{eqnarray*}
  KLD\big( g(t) \big| \! \big| f(t|s) \big)
   & \stackrel{s \rightarrow \infty}{\sim} &
      k + \frac{k}{2} \ln \pi + \ln \Gamma(k+1)  - \ln \Gamma(k/2+1) -
      k \ln s \nonumber \\
   & \stackrel{s \rightarrow \infty}{\sim} &
      c(k) - k \ln s.
\end{eqnarray*}
That is, the deviation for large $s$ of the length distribution of
the Waxman graph from the GER graph is a constant dependent on the
dimension of the space minus a term proportional to the log of $s$.
 

The likelihood function for a particular Waxman graph under the usual
independence assumption, given the lengths of the observed links
${\mathbf d}$, is
\begin{equation}
  \label{eq:lieklihood}
    {\cal L}( s \mid {\mathbf d}) 
   =  \prod_{(i,j) \in {\cal E}} f(d_{i,j} \mid s ) .
\end{equation}

We apply \autoref{eq:lieklihood} below, even though links in the
Waxman graph are only conditionally and not truly independent given the distances. 
This can be intuited by considering triangles
(edges connecting sets of three nodes). The fact that the underlying distances are a
metric forces the three distances to satisfy the triangle inequality
and thus the three must be correlated. 
Fortunately the correlation is largely local. For instance,
if we consider two distinct node pairs $(i_1,j_1)$ and $(i_2,j_2)$,
then the existence of edges $(i_1,j_1)$ and $(i_2,j_2)$ is independent
(given no extra information). Thus correlations are largely mediated
through common nodes.
As a result we should expect weaker correlations as the number of
nodes $n$ grows (given a fixed number of edges) and similarly we
expect the correlations to increase as the average node degree
increases.  We demonstrate how this affects various estimators below.

\section{Estimation Techniques}

This section describes various approaches to the parameter estimation
of Waxman random graphs.

\subsection{Log-linear regression}

The first estimation method applied to this problem was introduced by
\citet{Lakhina:2002:GLI:637201.637240}.  They noted that
\begin{equation}
  \label{eq:deterrence_fn}
  \frac{ f(d \mid s) }{ g(d) } = c \exp(- s d),
\end{equation}
where we can see from our calculations that $c = 1/\tilde{G}$.  If we
could form this ratio, we might estimate $s$ by log-linear regression
against $d$. The observed lengths of links in the graph yield implicit measurement 
of the numerator in \autoref{eq:deterrence_fn} and we could obtain estimates of $g(d)$ either:
\begin{enumerate}

\item {\em analytically:} use the shape of the region to compute
  $g(d)$; or

\item {\em empirically:} use the distances between all of the nodes
  (not just the linked nodes) to estimate $g(d)$.

\end{enumerate}
The former has the advantage of being fast. The latter exploits
the data itself in the case that the region is not regular, but
computationally it is $O(n^2)$.  The latter approach was used by
\citet{Lakhina:2002:GLI:637201.637240}, but we shall evaluate both.

We also need to estimate $f(d \mid s)$:
\citet{Lakhina:2002:GLI:637201.637240} did so by forming a
histogram. We shall use this approach, so the estimator
proceeds by counting the number of edges in each length bin, and
dividing by the expected number in that bin absent the distance
deterrence function.

Once we have computed $\hat{s}$, we estimate $q$ by inverting
\autoref{eq:average_edges} to get
\begin{equation}
  \label{eq:q_est}
  \hat{q} = \frac{2 e}{n(n-1) \tilde{G}(\hat{s})},
\end{equation}
where $n$ and $e$ are the observed number of nodes and edges
respectively.  The decoupling of estimation of $s$ and $q$ makes this
sequential estimation possible, and this will be exploited in other
estimators described below.

The time complexity of the above algorithm based on the analytical formulation of $g(d)$
is $O(e)$  but if we were to estimate $g(d)$ using the
distances between all nodes then it would be $O(n^2)$. 
Our evaluation focuses on the time to perform the
regression, \ie we do not include the time to construct the distances
and empirical histograms as these are problem dependent. For
instance, Euclidean distance calculations are fast, but in one case we
compute distances on the globe (including corrections for its
non-spherical nature) and these calculations can be considerably slower. 

\subsection{Generalised Linear Model (GLM)}

\citet{davis14:_spatial} use a GLM to estimate the parameters of their
model, ours is slightly different. In order to make the difference in usage explicit we 
include a description of the framework taken from  \citep[p.~591]{Casella2002}.
 A GLM consists of three components: 
\begin{description}
\item[A random component:] The response variables $Y_1, Y_2, \ldots,
  Y_n$. These are assumed to be independent (but not identically
  distributed) random variables from a specified exponential family.

\item[A systematic component: ] a linear function of a vector of
  predictor variables
  \begin{equation}
    \beta_0 + \sum^k_{m= 1} \beta_m \mathbf{X}^{(m)},            
    \label{eq:glm_1}
  \end{equation}
  where $\mathbf{X}^{(m)}$ is the $m$th predictor variable.
  
\item[A link function:]  $r(\mu)$  such that
  \begin{equation}
    r(\mu_i) =  \beta_0 + \sum^k_{m= 1} \beta_m \mathbf{X_i}^{(m)},
    \label{eq:glm_2}
  \end{equation}
  where $\mu_i = \E{Y_i}$, and $\mathbf{X_i}$ is the predictor
  vector for the $i$th response. 
\end{description}

Estimation of the $\beta_m$ is usually by maximum likelihood, \ie we
write the log likelihood, differentiate with respect to $\beta_m$ and
set to zero, giving $k+1$ non-linear equations which are solved
numerically. Here we use Matlab's {\tt glmfit} to perform the task.


We index the variables by their edge label $(i,j)$ and the response
variables $Y_{(i,j)}$ are indicators of edges so $\mathbf{Y}$ is the
adjacency matrix of the graph, \ie
\begin{equation}
  \label{eq:indicator_edge}
  Y_{(i,j)} = \left\{
    \begin{array}{ll}
      1, & \mbox{ if } (i,j) \in {\cal E}, \\
      0, & \mbox{ otherwise.}
    \end{array}
    \right.
\end{equation}
These are modelled as  Bernoulli random variables, and hence $\mu_{(i,j)} =
p(d_{i,j})$, \ie the probability that edge exists given that the distance
between nodes $i$ and $j$ is $d_{i,j}$.  The predictors
are the distances, \ie $X_{(i,j)} = d_{i,j}$ and $k=1$.

The probabilities are assumed to follow \autoref{eq:waxman_prob},
so the natural link function is $\log$. In this case we get
\begin{equation}
  \log( p_{(i,j)} ) =  \beta_0 + \beta_1 d_{i,j},
\end{equation}
and  make the natural identification that $\beta_0 = \log(q)$ and
$\beta_1 = - s$.

The GLM of \citet{davis14:_spatial} is almost identical but for the
use of the link function $r(p_{(i,j)}) = \log(-\log(1 - p_{(i,j)}))$.
 
When the response random variables are Bernoulli, it is more conventional to use a
logit link function. The advantage of this is that it maps
the entire range of linear combinations of the systematic component
into the interval $(0,1)$, and thus produces legitimate
probabilities. 

Using $\log$ as the link function for modelling Waxman graphs  
means there is an implicit bound on the parameters 
$\beta_m$ such that $\beta_0 + \beta_1 d_{i,j}$ must be
positive. As long as the estimator returns positive parameters, this
is assured, but in cases such as small $s$ it is possible
to have examples where natural estimates lie on the boundary. 
Similar issues arise for our estimator, so we defer further discussion 
until section \ref{existence}.

One potential problem with fitting a GLM is that if there is a
cutoff $c$ such that for all $d_{ij} \leq c$ the edge exists and for
all $d_{ij} > c$ the edge is absent, then the GLM may not find finite
coefficients that best fit the data and estimates will have
large errors. This is unlikely to happen when estimating the parameters of 
large networks since the probability as a function of $d$ decays exponentially,
so there is no clean cutoff. However, in small networks there may
only be a few edges, and it is possible that the data exhibit
such a cutoff by chance. In cases with $n \sim 10$ we observe very
poor performance for GLM fits as a result of this phenomena, but it has
negligible impact on the accuracy of GLM fits for larger networks.

The GLM method uses all possible edges: those that exists and those
that do not. Hence the number of response variables in the problem is
$O(n^2)$, and so the computational time complexity is also $O(n^2)$.

\subsection{Sufficient Statistics}

An obvious question arises as to what is a set
of {\em minimal sufficient statistics} for use in the estimation of Waxman random graph parameters. 
To answer this we apply the following theorem:

\begin{thm}[\citet{Casella2002}, Theorem 6.2.13, p.281] 
  \label{thm:min_suf}
  \mbox{ } \\
  Let $f({\mathbf x} \mid {\mathbf \theta})$ be the PMF or PDF of a
  sample $\mathbf X$. Suppose there exists a function $T({\mathbf x})$
  such that, for every two sample points $\mathbf x$ and $\mathbf y$,
  the ratio $f({\mathbf x} \mid {\mathbf \theta}) / f({\mathbf y} \mid
  {\mathbf \theta})$ is a constant function of $\mathbf \theta$ if and
  only if $T({\mathbf x}) = T({\mathbf y})$. Then $T({\mathbf X})$ is
  a minimal sufficient statistic for $\mathbf \theta$.
\end{thm}

Take the sample to be the set of edge lengths ${\mathbf d} = (d_1,
\ldots, d_e)$ where $d_k = d_{(i,j)}$ for $(i,j) \in {\cal E}$.
Under the independence assumption the PDF of a sample is
$\prod_{i=1}^e f\big(d_i | s, q\big)$, conditional on the number of
edges $e$, where $f$ is defined in \autoref{eqn:link_len_dist}, and
$e$ is binomially distributed $B^N_p(e)$ where $N = n(n-1)/2$ and $p =
q \tilde{G}(s)$.

Consider two samples: ${\mathbf x}=(x_1,\ldots,x_{e_1})$ and ${\mathbf
  y}=(y_1, \ldots, y_{e_2})$, then the ratio in \autoref{thm:min_suf}
is
\begin{eqnarray}
  \frac{  B^N_p(e_1) \prod_{i=1}^{e_1} f\big( x_i | s, q \big)}
       {  B^N_p(e_2) \prod_{i=1}^{e_2} f\big( y_i | s, q \big)}
 & = & \frac{   B^N_p(e_1) \prod_{i=1}^{e_1} g(x_i) e^{-s x_i} /\tilde{G}(s)}
            {   B^N_p(e_2) \prod_{i=1}^{e_2} g(y_i) e^{-s y_i} /\tilde{G}(s)} \nonumber \\ 
 & = & \frac{ \prod_{i=1}^{e_1}  g(x_i)}
            { \prod_{i=1}^{e_2}  g(y_i)}  \times
       \frac{ \prod_{i=1}^{e_1}  e^{-s x_i} }
            { \prod_{i=1}^{e_2}  e^{-s y_i} } \times
        \frac{ \tilde{G}(s)^{e_2} B^N_p(e_1)}{ \tilde{G}(s)^{e_1} B^N_p(e_2)}
          \nonumber \\ 
 & = &   m({\mathbf x},{\mathbf y}) e^{-s (e_1 \bar{d}_x - e_2 \bar{d}_y)}
       \tilde{G}(s)^{e_2 - e_1}
       \frac{B^N_p(e_1)}{B^N_p(e_2)},
\end{eqnarray}
where $\bar{d}_x$ and $\bar{d}_y$ are the averages of the distances in
datasets ${\mathbf x}$ and $\mathbf y$ respectively, and $m({\mathbf
  x},{\mathbf y})$ is a function only of $g(\cdot)$ of ${\mathbf x}$
and ${\mathbf y}$, independent of the parameters $(q,s)$.

The ratio above depends on the parameters $s$ and $q$ only through the
statistics $\bar{d}$ and $e$. Hence if $e_1=e_2$ and $\bar{d}_x =
\bar{d}_y$, the ratio is a constant function of $(q,s)$. On the other
hand, if either $e_1 \neq e_2$ or $\bar{d}_x \neq \bar{d}_y$, then the
ratio is a non-constant function of the parameters. Hence the
conditions of the theorem are satisfied and $(e,\bar{d})$ forms a
minimal sufficient set of statistics.  Notably, these statistics can
be constructed almost trivially in $O(e)$ time, and we will use that
fact to construct a MLE.


\subsection{Maximum Likelihood Estimator}

We consider the log-likelihood function ${\ell}(s \mid {\mathbf d}
) = \ln {\cal L}(s \mid {\mathbf d})$ and note from
\autoref{eq:lieklihood} that
\begin{eqnarray*}
 {\ell}(s \mid {\mathbf d} ) 
   & = &  \sum_{{(i,j)} \in {\cal E}} \ln  f(d_{(i,j)} \mid s )  \\
   & = &  \sum_{{(i,j)} \in {\cal E}} \left[ \ln g(d_{(i,j)}) - s d_{(i,j)} - \ln \tilde{G}(s) \right] \\
   & = &  - e \ln \tilde{G}(s) + \sum_{{(i,j)} \in {\cal E}} \ln g(d_{(i,j)}) -  s \sum_{{(i,j)} \in {\cal E}} d_{(i,j)}.
\end{eqnarray*}
Our goal is to find $s$ such that the partial derivative of
${\ell}(s \mid {\mathbf d} )$ w.r.t. $s$ is zero, \ie
\begin{equation}
  \frac{\partial}{\partial s} {\ell}(s \mid {\mathbf d}) 
    =  - e\, \frac{\tilde{G}'(s)}{\tilde{G}(s)} - \sum_{{(i,j)} \in {\cal E}} d_{(i,j)}
    = 0 ,
    \label{eq:likelihood_d}
\end{equation}
so we need to find $s$ such that
\begin{equation}
 \frac{\tilde{G}'(s)}{\tilde{G}(s)} = - \frac{1}{e}  \sum_{{(i,j)} \in {\cal E}} d_{(i,j)}
     = - \bar{d} , 
 \label{eqn:fund_est}
\end{equation}
where $\bar{d}$ is the observed mean link length. 

From (\ref{eqn:av_link_dist}) we know $-\tilde{G}'/\tilde{G}$ is the
expected length of line segments on the Waxman graph, so the MLE of
$s$ is also the moment-matching estimator.

\subsection{Existence and uniqueness of the MLE}
\label{existence}
When will a unique solution to \autoref{eqn:fund_est} exist?  We know
from the definitions of $\tilde{G}(s)$ and its derivative that
$\tilde{G}(s)>0$, and $h(s) = -\tilde{G}'(s)/ \tilde{G}(s)$ is a
continuous, positive function for all $s$, so if $h(\cdot)$ is
monotonic there will be at most one solution to \autoref{eqn:fund_est}
for any given $\bar{d}$.

Consider the derivative 
\begin{equation}
  \label{eq:derivative_h}
  \frac{dh}{ds} =  \frac{\tilde{G}'(s)^2 - \tilde{G}''(s) \tilde{G}(s)}{\tilde{G}(s)^2}.
\end{equation}
Now the numerator is positive, and from Schwarz's Inequality 
\begin{eqnarray*}
  \tilde{G}''(s)\tilde{G}(s) & = & \int \left|\sqrt{t^2 g(t) e^{-s t} } \right|^2\,dt \cdot
                     \int\left| \sqrt{g(t) e^{-s t}}\right|^2\,dt  \\ 
           & > & \left|\int t g(t) e^{-s t} \,dt\right|^2 \\
           & = & \tilde{G}'(s)^2 ,
\end{eqnarray*}
and so the denominator is negative, and hence the function $h(s)$ is
monotonically decreasing with $s$.


When $s=0$, the Laplace transform $\tilde{G}(0)=1$ and hence
$h(0)=\int t g(t) dt = \bar{g}$, the average distance between pairs
of nodes. When we remove longer links preferentially, the average link
distance must decrease, so it is intuitive that $h(s)$ is a decreasing
function, starting at $h(0) = \bar{g} = d_{max}$, which is the maximum expected
edge distance over all possible parameters $s$.

In the limit as $s \rightarrow \infty$ we already know from
\autoref{eq:large_s_d} that $h(s) = \E{d} \sim k / s$, and so we know
that in the limit as $s \rightarrow \infty$ that $h(s) \rightarrow 0$,
so for any measured $\bar{d} \in (0, d_{max}]$ there will be a unique
solution to \autoref{eqn:fund_est}.

Unfortunately, it is possible for the sample mean of the edge
distances $\bar{d} > d_{max}$, \ie for a particular graph to have an
unlikely preponderance of longer links. In this case
\autoref{eqn:fund_est} has no solution for $s \in [0, \infty)$, which
seems to create problems. However, the obvious interpretation of
$\bar{d} > d_{max}$ is that then there is no evidence that long links
have been preferentially filtered from the graph, and so it is natural
in this case to assume the model should be the GER random graph, \ie
$s=0$.

In formal terms, the MLE satisfies some but not all of the properties required for
it to be {\em consistent}. Standard asymptotic theory for
MLEs requires the condition that the true parameter value lies away
from the boundary to form consistent estimates that converge to the
true value as the amount of data increases.

For more formal consideration of this boundary case, we could
introduce a strict hypothesis test for GER vs Waxman graphs, but for
the moment we take $s=0$ in these cases. We will see that in these
cases that this leads to a slight bias in estimates for graphs with
small $s$. The bias leads to a variance-bias tradeoff, leading to a
RMS error in the estimator that can be lower than the Cram\'{e}r-Rao
(CR) bound (for unbiased estimators).

In the case where the parameter $s$ does lie away from the boundary,
ideal MLEs have a number of useful properties:
\begin{itemize}

\item {\em Consistency:} a sequence of MLEs converges in probability
  to the value being estimated as the amount of data increase: in
  particular the estimates are asymptotically unbiased.

\item {\em Efficiency:} \ie it achieves the CR lower bound when the
  sample size tends to $\infty$, \ie no consistent estimator has lower
  asymptotic mean squared error.

\end{itemize}
Another way to state these results is that MLEs are {\em
  asymptotically normal}, \ie as the sample size increases, the
distribution of the MLE tends to the Gaussian distribution with mean
given by the true parameter, and covariance matrix equal to the
inverse of the Fisher information matrix.

However, we know that our MLE does have a potential problem at the
boundary, and that the independence assumption is an approximation,
and so it is interesting to consider whether the ideal bound is
actually achieved. In the next section we derive the CR lower bound
against which we make comparisons.

\subsection{Bounds}

The CR lower bound is the minimum variance of any unbiased estimator,
under certain mild conditions. In the case of a scalar parameter
such as $s$, the bound is given by the inverse of the Fisher
information $I(s)$ defined by
\begin{equation}
  \label{eq:cr_bound}
  I(s)
   = - \E{ \left( \frac{\partial^2 \ell(d; s)}{\partial s^2}  \right) }.
\end{equation}
Here we have already computed the first derivative of the likelihood
in \autoref{eq:likelihood_d}, and its second derivative is
\begin{equation}
   \frac{\partial^2 \ell}{\partial s^2} 
   =  e\, \frac{dh}{ds},
   \label{eq:likelihood_dd}  
\end{equation}
where we calculated $dh/ds$ in \autoref{eq:derivative_h}. Taking the
expected value, and noting 
\autoref{eq:average_edges}, we get
\begin{equation}
 I(s)
    = \frac{n (n-1) q}{2} \left[ \frac{\tilde{G}''(s) \tilde{G}(s) - \tilde{G}'(s)^2}{\tilde{G}(s)^2} \right].
\end{equation}
The net effect of this result is that an estimator that achieves the
CR bound would have variance $O(1/n^2)$. 

However, that assumes that growth of the number of edges is $O(n^2)$,
but a more reasonable model might increase the size of the network
while keeping the average node degree $\bar{k}$ constant, in which
case \autoref{eq:node_degree} implies that the variance will be $O(1 /
n)$. Examples of this can be seen in the results.  

Even if the network is not sparse, if we randomly sample $O(n)$ edges,
then the result above implied that the  variance will be $O(1/n)$,
with the benefit that (for some estimators) the cost of computation
will now be $O(n)$. 

The mean squared error of an estimator is composed of the variance
plus the bias-squared, so it is possible for a biased estimator to
improve on the bound (as we shall see in special cases below), but in
general we measure estimator efficiency (with respect to accuracy as a
function of the amount of data) with respect to this bound.


\subsection{MLE of $\hat{q}$} 

As noted above, once we know $\hat{s}$, we can use this to estimate
$\hat{q}$ using \autoref{eq:q_est}. If we know the exact value of $s$
then $q$ acts as a random filter of links, so estimation of $q$
is equivalent to estimating the parameter of a Bernoulli process. 
Thus \autoref{eq:q_est} would be the MLE of $q$ given the true
value of $s$. The question is whether it is still a MLE when we derive
it from the estimated value of $\hat{s}$.

To answer this question we draw on the property of functional
invariance of MLE \ie if parameters are related through a
transformation $a = g(s)$, then the MLE of $a$ is $\hat{a} = g(
\hat{s} )$.  This property allows  the conversion back to the 
original Waxman parameter $\alpha = 1/ s L$, for instance.
More importantly it means that $\hat{q} = 2e /
n(n-1) \tilde{G}(\hat{s})$ is the MLE of $q$, given a MLE of
$\hat{s}$.

Once again, the property above applies for the ideal MLE. We test the
estimator $\hat{q}$, and we show that its properties are close to
those of the ideal MLE, except when $s$ is small.

\subsection{Numerical calculation of the MLE}

MLEs are often computed directly using techniques such as Newton's
method. In this problem, each iteration would require computation of a
Laplace transform\footnote{The Laplace transforms are calculated
  numerically using Matlab's {\tt integral} function, which performs
  adaptive numerical quadrature.  }, and its derivative, and so will
be somewhat slower than many MLEs. However, we can improve this
performance in several ways. First, we rearrange
\autoref{eqn:fund_est} in the form
\begin{equation}
  \label{eq:fast_mle}
  \bar{d} \tilde{G}(s)  + \tilde{G}'(s)
  = {\cal LT}\big[ (\bar{d} - t) g(t); s \big] 
  = 0,
\end{equation}
which reduces the number of transforms that need to be calculated.

The estimation algorithm is a 1D search using Matlab's {\tt fzero}
function. No doubt faster implementations are possible using
Newton-Raphson or other approaches, but we found this to converge in
only a few steps, given reasonable starting bounds, which we choose to
be $[0, k/\bar{d}]$ (where the upper bound is given by the asymptotic
form of the average distance \autoref{eq:large_s_d}).

In the particular case that $\bar{d} > d_{max}$, \ie the sample mean
of the edge distances is greater than the largest expected edge
distance, then we take $\hat{s} = 0$. 

The speed of calculation can be improved if we precalculate $g(t)$ at
a set of fixed grid points, and reuse this in all of the Laplace
transform calculations, thus avoiding frequent recomputation of the
density. We create a second estimate, the MLE-N (Numerical), in which
we calculate the inverse CDF at uniformly spaced probabilities to
minimise errors arising from the fixed grid.  In the test below we use
a grid of 1000 points. We will see that this approach allows faster
computation at the cost of a small increase in error for large
$s$. The tradeoff could be adjusted by adapting the number of points
used.

The other reason for using points spaced uniformly in the probability
space is that it is relatively easy to estimate an inverse CDF from
data. In particular we could, if available, take the complete set of
distances between all nodes and derive from this an empirical estimate
of the inverse CDF to use in the MLE calculations. We refer to this
method as the MLE-E (Empirical). This method avoids the problem of
estimating the region size and shape, and the assumption that nodes
are uniformly distributed over the space.
 
\section{Tests} 
 
In this section we examine the performance of the Waxman graph
parameter estimators on simulated data. There are three key parameters
for a Waxman graph, $s$, $q$ and $n$. Instead of using these we choose
$s$, $n$, and the average node degree $\bar{k}$ which implicitly
determine $q$. We primarily focus on fairly sparse networks
($\bar{k}=3$) because sparsity is the common case for many real
networks, but we do examine the effect of changes in these parameters
(note there are some combinations of these parameters for which $q$
cannot exist).

The accuracy of the methods increase as the graphs become larger so we choose parameter values that result in moderately sized (1000 node) graphs.  This ensures that errors are of a magnitude that allows us to compare and assess them. However, we do also examine the relationship between graph size and accuracy of parameter estimates. 

Except where otherwise noted we simulate 1000 Waxman random graphs for each
parameter setting using Matlab, and form estimates of bias and RMS
(Root-Mean-Square) error for each estimation method. 
All unannotated results are for the Waxman graph on the unit square
with the Euclidean distance metric. 


\subsection{Comparisons}

\begin{figure}[tbp]
  \begin{center}

    \includegraphics[width=\figurewidthA]{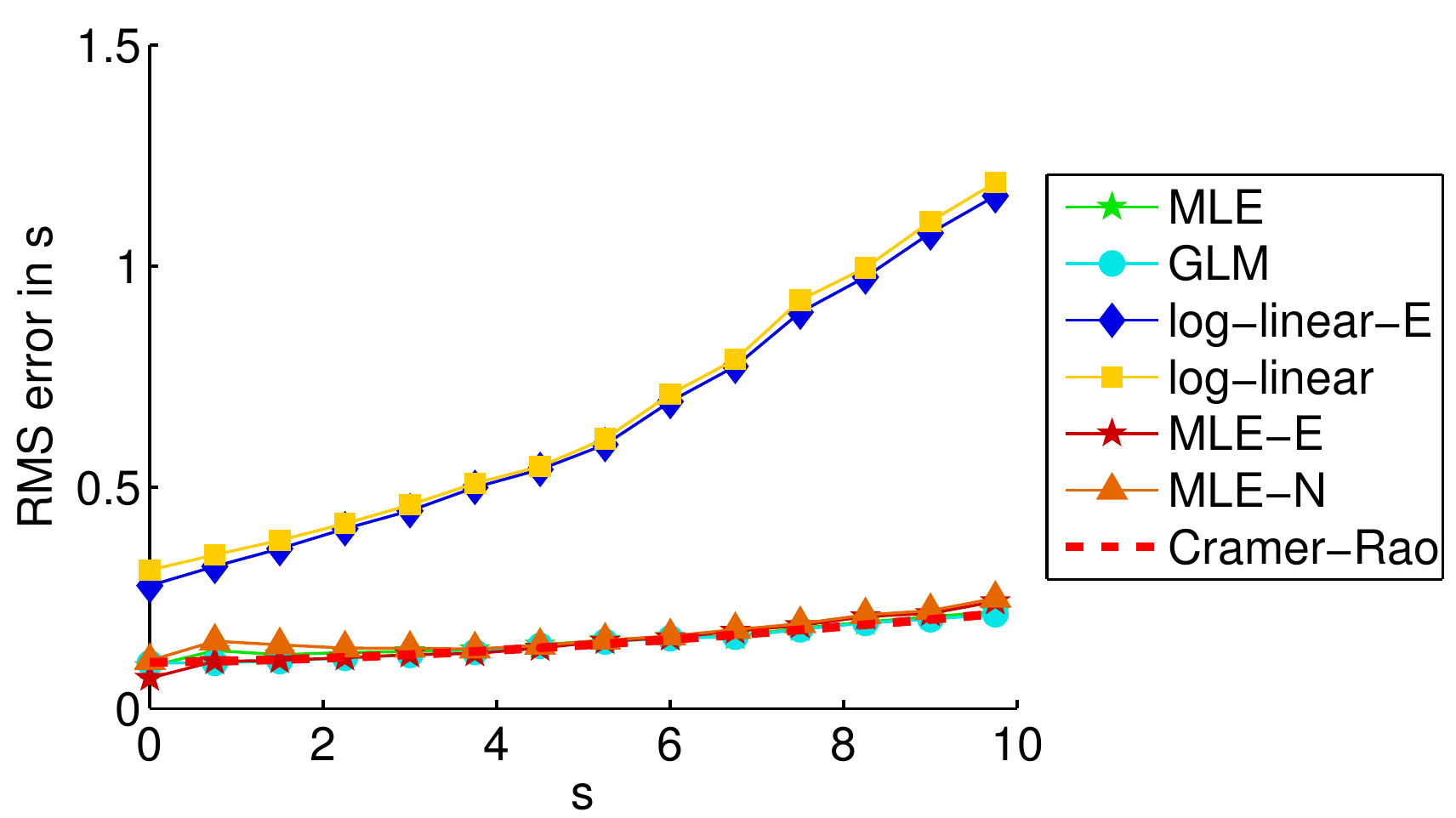}
      \caption{RMS errors as a function of $s$ ($\bar{k}=3$,
        $n=1000$). See \autoref{fig:rms_s_z} for detail.}
    \label{fig:rms_s}
  \end{center}
\end{figure}

We start with a basic test of the estimation methods. We consider the
RMS (Root Mean Squared) errors of the methods with respect to the CR
lower bound as a function of the parameter $s$. \autoref{fig:rms_s}
shows these errors. The most obvious conclusions are that the
log-linear regression is the least accurate and that the GLM and MLE
estimators are all roughly equivalent in accuracy. The minor empirical
variants of the log-linear and MLE estimators apparently have only
small affects on accuracy.

\begin{figure}[tbp]
  \begin{center}
 
    \includegraphics[width=\figurewidthA]{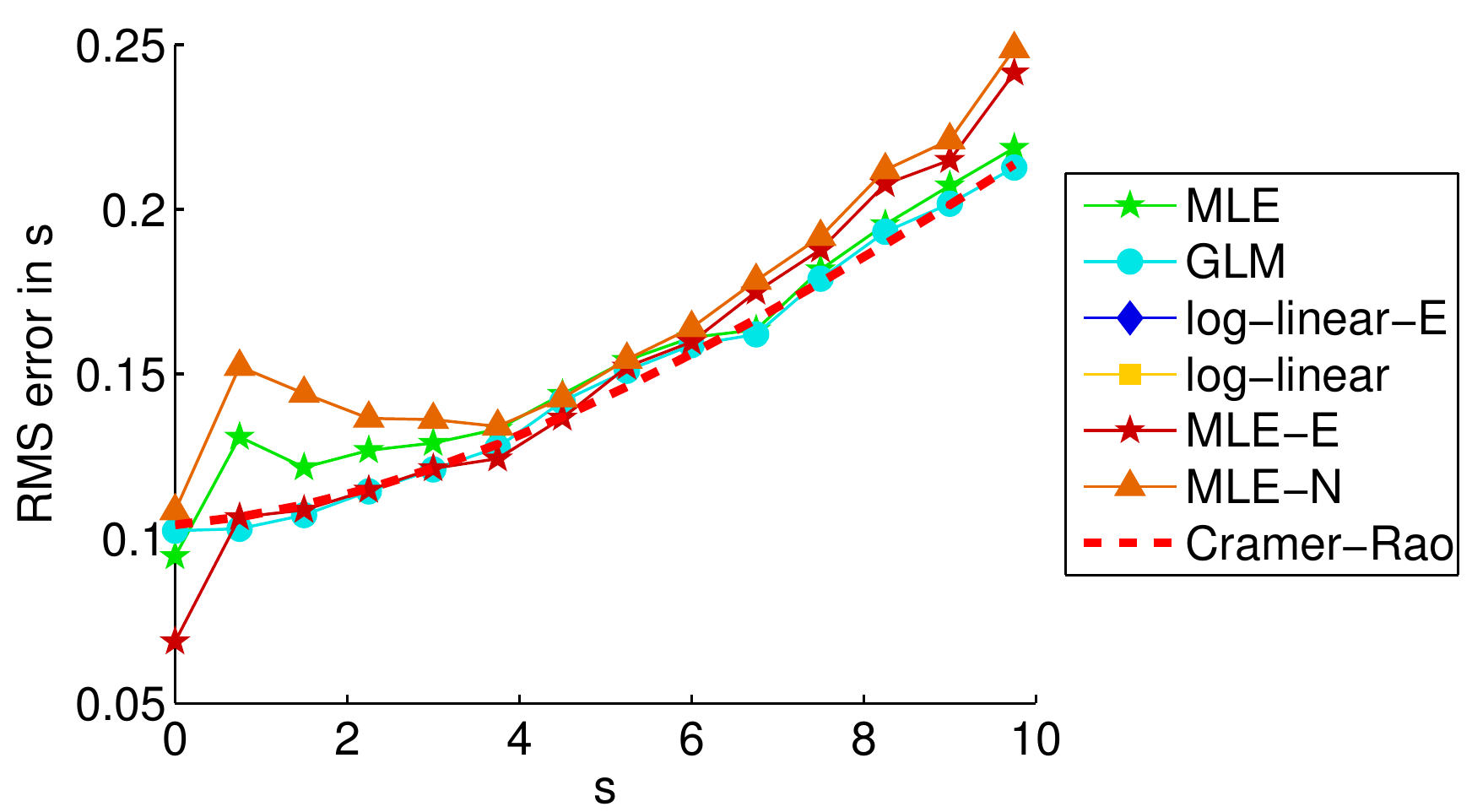}
      \caption{RMS errors as a function of $s$ ($\bar{k}=3$, $n=1000$)
        detail.}
    \label{fig:rms_s_z}
  \end{center}
\end{figure}

However, when we examine a more detailed view in \autoref{fig:rms_s_z}
some differences are apparent. Over the entire range the GLM tracks the CR
bound, but for small values of $s$ there is some advantage in using the GLM
and MLE-E. For very small $s$, the MLE approaches can actually beat
the CR bound. This can be explained by the constraint $s
\geq 0$. This constraint leads to a small bias in the MLE estimators
for small $s$. When $s$ is around 1, the bias leads to an increase in
the RMS errors. When $s$ is very small the estimators have slightly
more information than assumed in our naive CR bound, and this can
reduce the variance of the estimator \citep{Gorman90CCRB}.

For large values of $s$ the MLE-E and MLE-N show a small increase in
error compared to the exact MLE, because for large $s$ the
grid size chosen is not small enough for accurate integrals to
be calculated. This effect could be mitigated by choosing a finer grid or by choosing a non-uniform grid 
with more detail around $t=0$ (albeit at  additional computational cost). We will leave the adaption of grid size
to future work. 

For small to moderate values of $s$, the MLE-E estimator is slightly
better than its exact cousin. This is valuable, but remember that the
MLE-E requires the information on all distances, not just those of the
existing edges, whereas the MLE and MLE-N estimators need only the
measurement of $\bar{d}$.  In fact for large graphs the $\bar{d}$ used
the exact MLE could be calculated from a sample of edges, so it is
possible to achieve sub-$O(e)$ performance for such graphs.  We have
examined a number of parameter settings and have observed the same
trends in the results and below we examine accuracy as a function of
other parameter values.

\begin{figure}[tbp]
  \begin{center}
    \includegraphics[width=\figurewidthA]{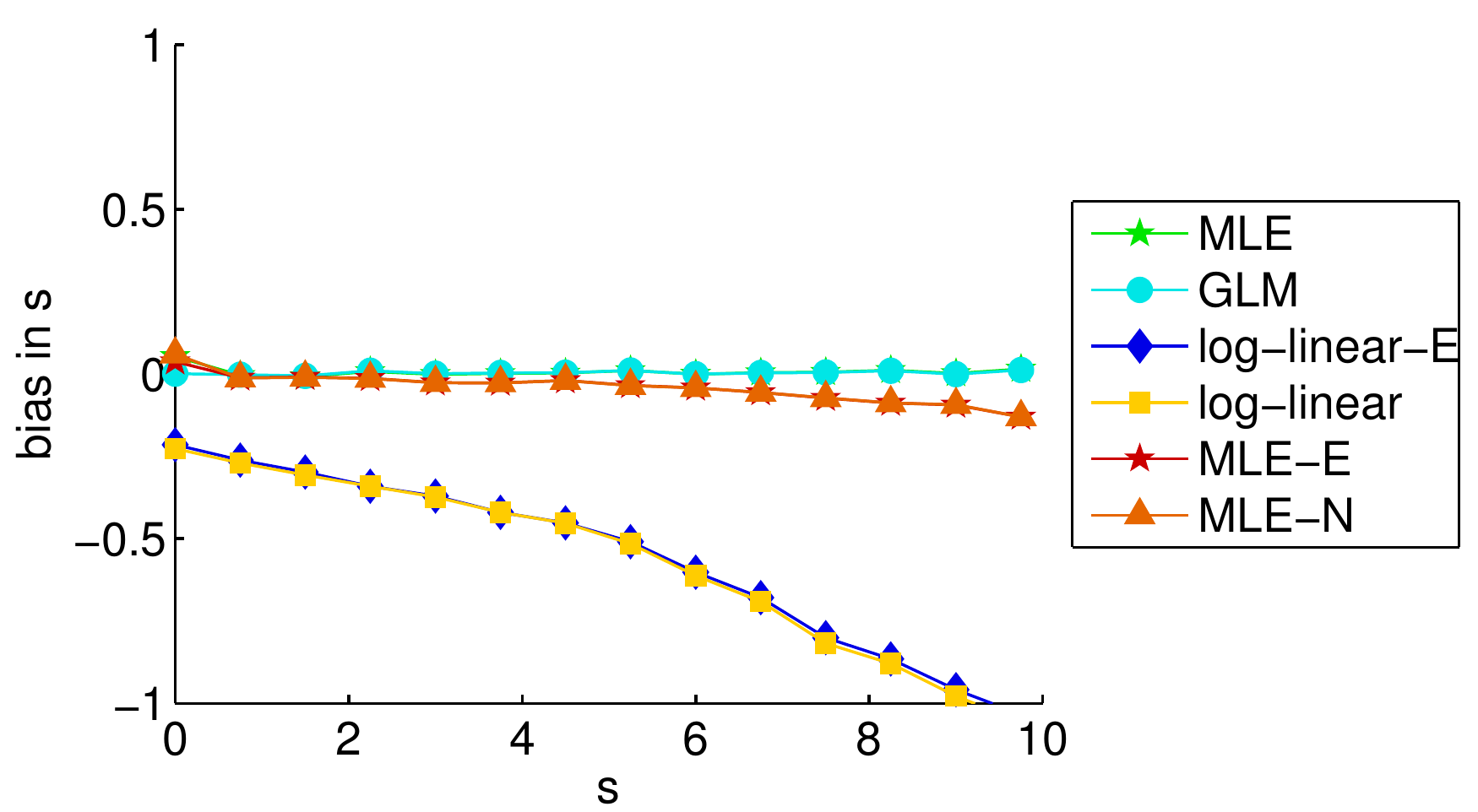}
    \caption{Bias in $\hat{s}$ for the different methods ($\bar{k}=3$,
      $n=1000$).}
    \label{fig:bias}
  \end{center}
\end{figure}

The major source of error can be see in \autoref{fig:bias}, which
shows the bias in the various approaches. We see immediately that the
log-linear approach suffers from significant bias which increases with
$s$ and remains the major source of error in the method. The
GLM is almost unbiased except for very small networks, where the
technique breaks down completely, resulting in very large errors.  On
the other hand, we see that the MLE is almost unbiased even for small
networks, except for small $s$ where there is a small positive bias
that explains the increase in RMS errors in that domain.

\begin{figure}[tbp]
  \begin{center}
    \includegraphics[width=\figurewidthA]{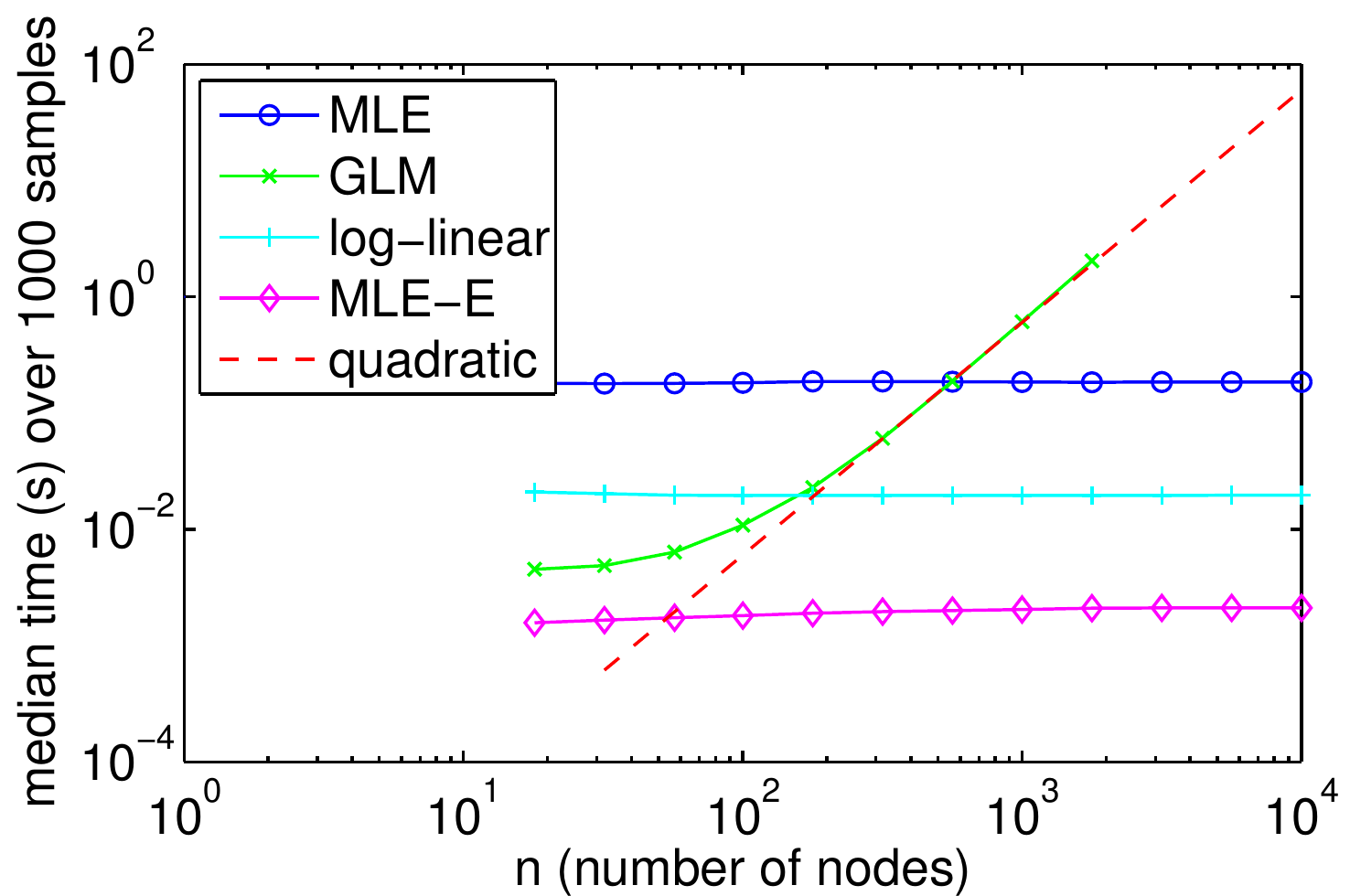}
    \caption{Computation times ($s=4$, $\bar{k}=3$). Note that only
      the estimator time is included, not the time to form the matrix
      of distances needed in the GLM and log-linear methods.}
    \label{fig:time}
  \end{center}
\end{figure}

We also estimate computation times using Matlab's {\em tic/toc} timers
to estimate wall-clock time of execution, and take the median over our
1000 samples to provide a robust estimate of the typical computation
times. Compute times are shown in \autoref{fig:time}, where we can see
that the GLM takes quadratic time. The log-linear-E method (not shown)
is also quadratic, but with a much small constant than the GLM,
because a histogram must be formed from all of the distances.

The other methods appear to have constant time with respect to the
network size $n$. Constant time is an illusion however: these methods
are actually $O(e)$, which for the networks considered is also
$O(n)$. The constant time component is obviously dominant for the
network sizes considered -- we would expect to see the linearity only
if we examine very large networks.

Memory requirement for the algorithms are
\begin{itemize}
\item GLM: $O(n^2)$;
\item log-linear (and -E variant): $O(b)$ where $b$ is the number of
  histogram bins;
\item MLE: $O(1)$ (assuming the numerical quadrature is
  memory efficient); and 
\item MLE-E and MLE-N: $O(c)$ where $c$ is the number of grid points. 
\end{itemize}
Note that all but the GLM use constant memory as a function of network
size, as they are based on summary statistics that can be computed via
a streaming algorithm that does not need to hold data in memory.

There is little doubt that if accuracy is the prime consideration then
the GLM is the best approach, whereas if there is a need for speed we
can trade off a small amount of accuracy, and use the MLE or
MLE-E. However, the practical considerations mean the GLM is not viable for
large networks. Extrapolating the computation times for the GLM shown
in \autoref{fig:time}, we can see that computation would take in the
order of 100 seconds for a graph of 10,000 nodes (and 15,000 links),
or around 3 hours for a graph with 100,000 nodes. Similarly the GLM
requires a large amount of memory for larger networks, as compared
with the almost trivial amount required in the other algorithms. On
the other hand, the accuracy of the log-linear approaches is poor.
Consequently we shall consider only the MLE-based estimates in detail
in what follows.

\subsection{MLE in detail}

\begin{figure}[tbp]
  \begin{center}
    \begin{subfigure}[t]{\figurewidthB}
      \centering
      \includegraphics[width=\textwidth]{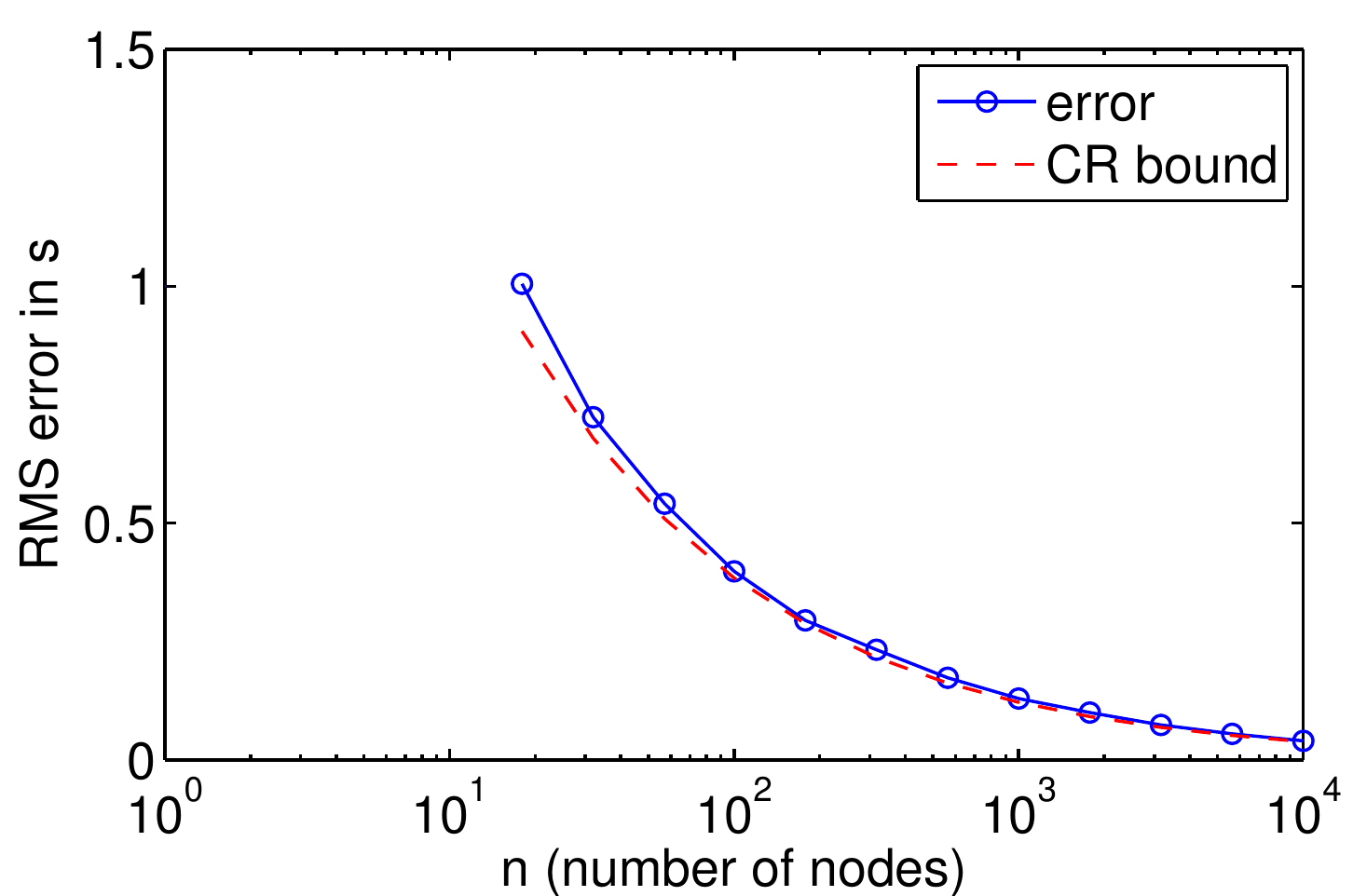}
      \caption{(a)}
    \end{subfigure}
    \hfil
    \begin{subfigure}[t]{\figurewidthB}
      \centering
      \includegraphics[width=\textwidth]{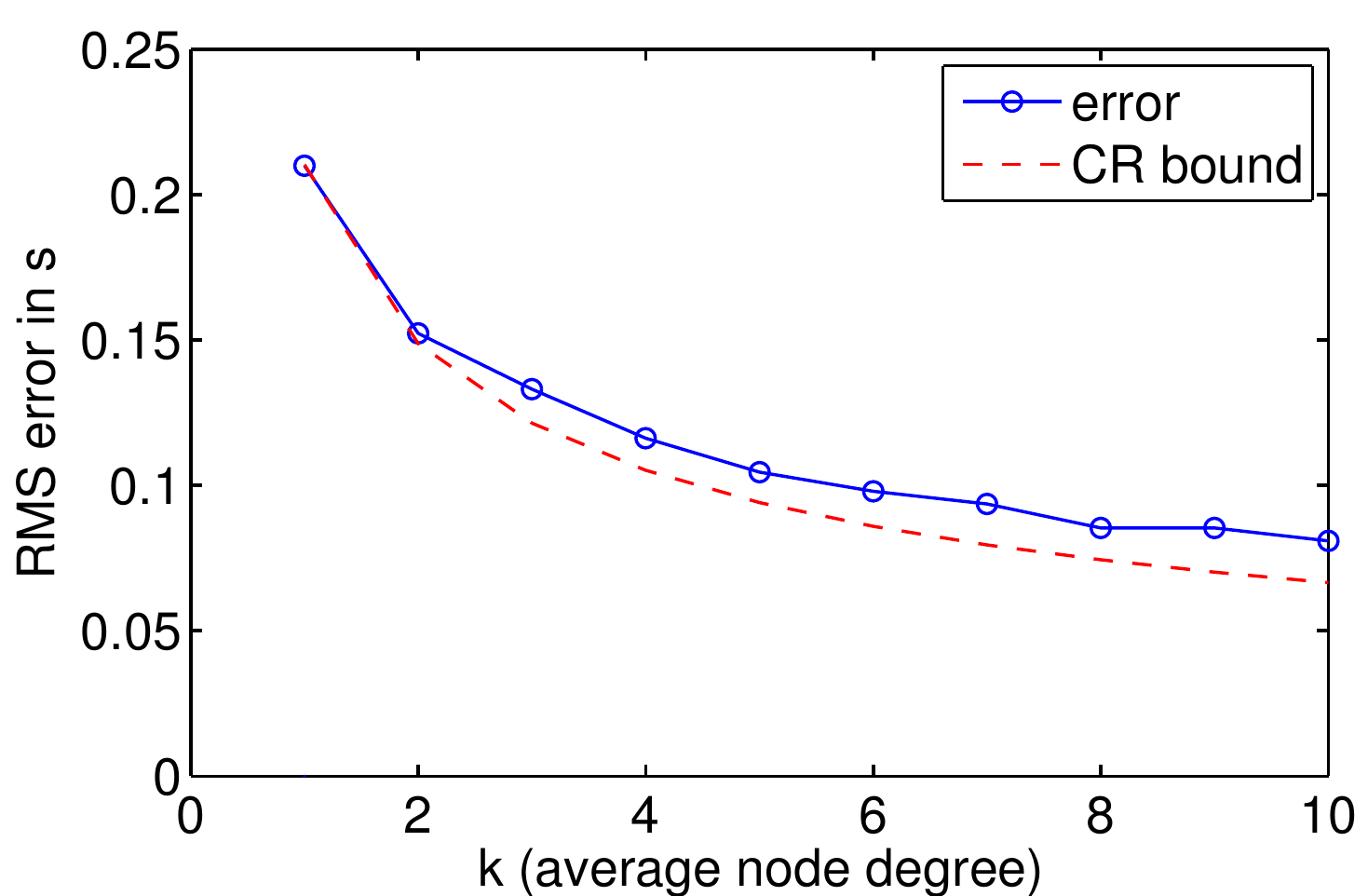}
      \caption{(b)}
    \end{subfigure}
    \caption{Accuracy of the MLE as a function of other parameters of the network
      (default parameters $n=1000$, $s=3.0$, and $\bar{k}=3$)      
      \label{fig:rms_n}
      \label{fig:rms_k}}
  \end{center}
\end{figure}

We now consider the accuracy of the MLE as both a function of the
network size $n$ in \autoref{fig:rms_n}(a) and of the average node
degree $\bar{k}$ in \autoref{fig:rms_k}(b). The accuracy of the method
is close to the CR bound for all parameter ranges, but suffers
slightly for small networks $n$, and for higher average node degree.

\begin{figure}[tbp]
  \begin{center} 
    \begin{subfigure}[t]{\figurewidthB}
      \centering
      \includegraphics[width=\textwidth]{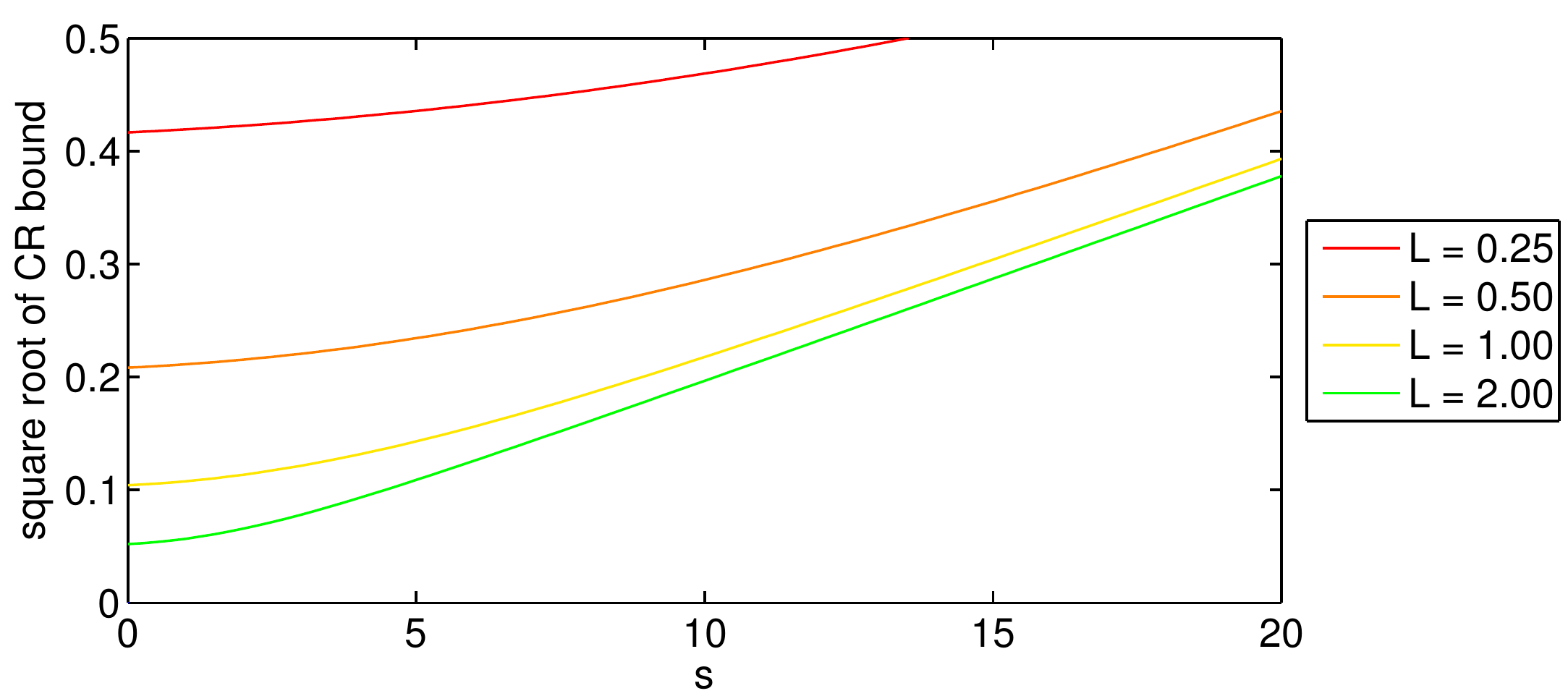}
      \caption{(a) Square region of size $L$.}
    \end{subfigure}
    \hfil
    \begin{subfigure}[t]{\figurewidthB}
      \centering
      \includegraphics[width=\textwidth]{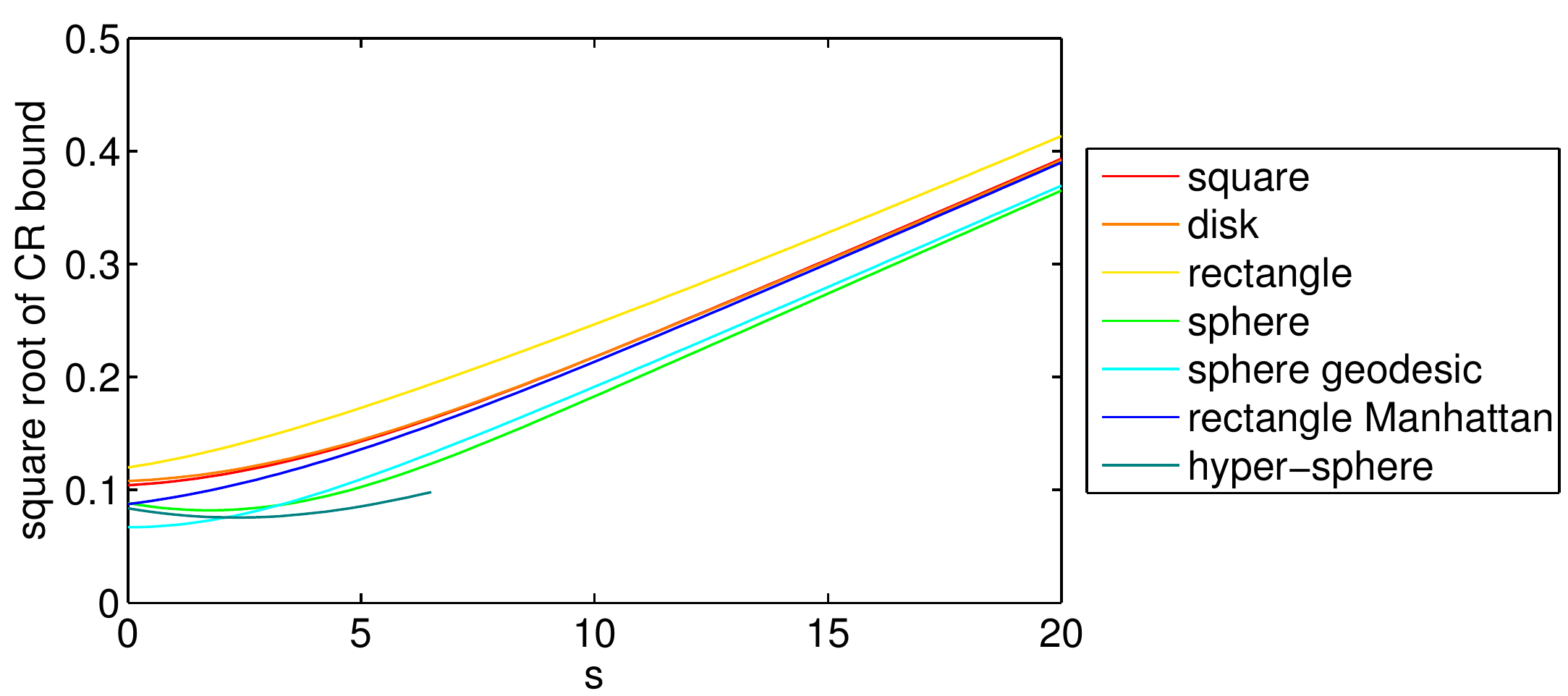}
      \caption{(b) Different region shapes.}
    \end{subfigure}
    \caption{CR bounds for different regions.
      \label{fig:rms_scale}
      \label{fig:rms_shape}
    }
    \label{fig:CR_regions}
  \end{center}
\end{figure}

All of the previous results consider accuracy of parameter estimates of the
Waxman graph over the unit square. Our approach naturally extends to a
variety of scales, regions, and distance
metrics. \autoref{fig:rms_scale}(a) shows the CR bounds for estimator
accuracy as a function of the length of the side of the square region,
and \autoref{fig:rms_shape}(b) shows these bounds for different region
shapes.

In the first case estimator accuracy is worse in smaller
regions. In such regions, all points are compacted closer together,
and there is a smaller range of possible length scales. Thus
we must estimate the exponential decay factor over a smaller
range of scales. The effect is more dramatic for smaller $s$ values,
where longer links are still likely. As $s \rightarrow \infty$ the CR
bounds converge because the length scale of typical links drops well
below the lengths of typical (potential) links. 

In the second case there are variations between regions but 
most have the same rough shape as a function of $s$. Notable
exceptions are sphere and hyper-sphere which are not monotonically
increasing functions of $s$ as the others are.

It is also notable that all of the methods have the same asymptotic
slope for large $s$ (though the function for the hyper-sphere is
truncated because larger $s$ is invalid given the fixed values of $n$
and $\bar{k}$ used in this plot). The overall conclusion is that the
shape of the region does have an effect (up to almost a factor of 2 in
estimator accuracy at some parameter values). However, the the shape
of the region is not the most important factor in determining the
accuracy of the estimator: the size of the graph, and its parameters
play a larger role. 

\subsection{Robustness}

\begin{figure}[tbp]
  \begin{center}
    \begin{subfigure}[t]{\figurewidthB}
      \centering
      \includegraphics[width=\textwidth]{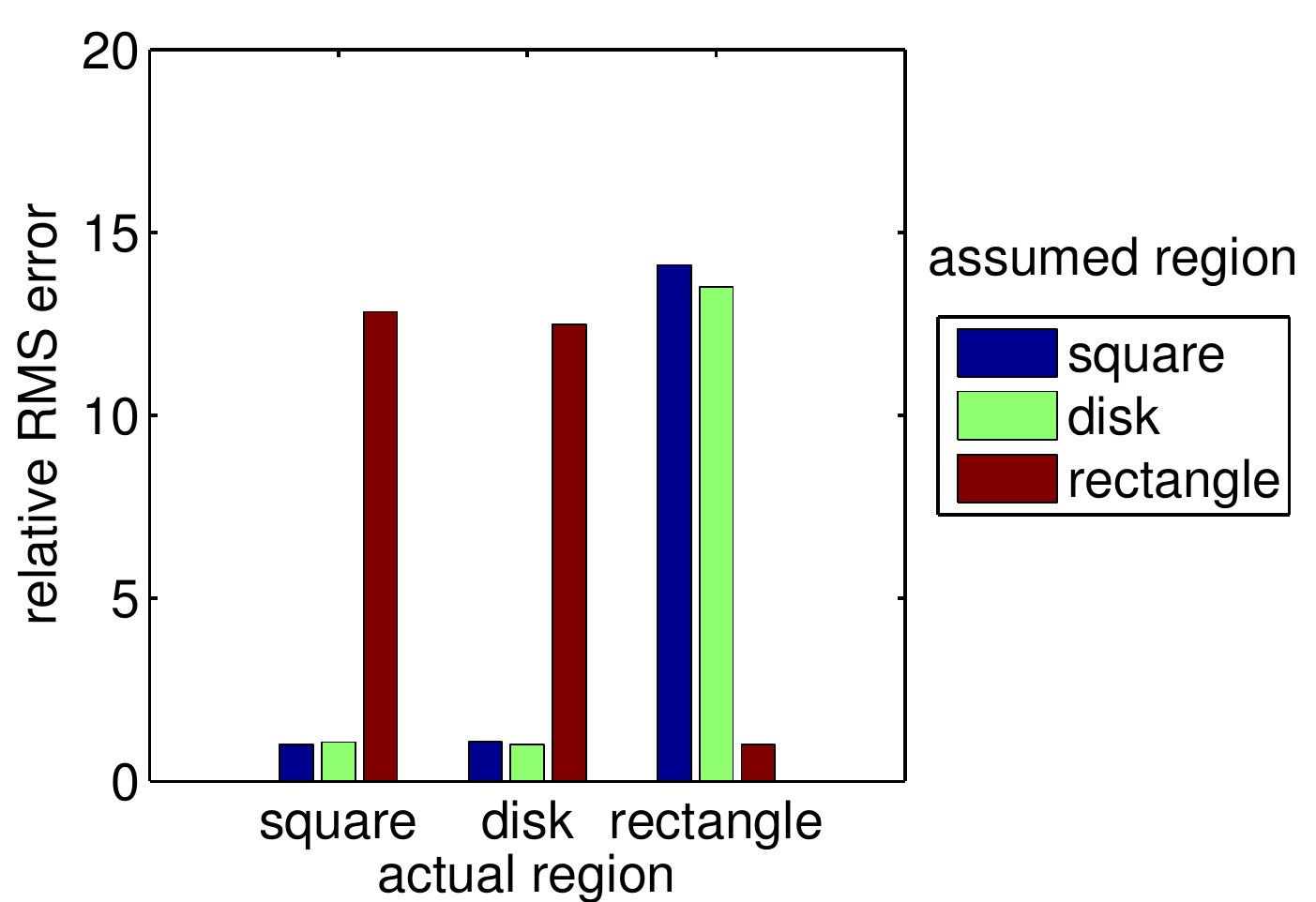}
      \caption{(a) $s = 3.0$.}
    \end{subfigure}
    \hfil
    \begin{subfigure}[t]{\figurewidthB}
      \centering
      \includegraphics[width=\textwidth]{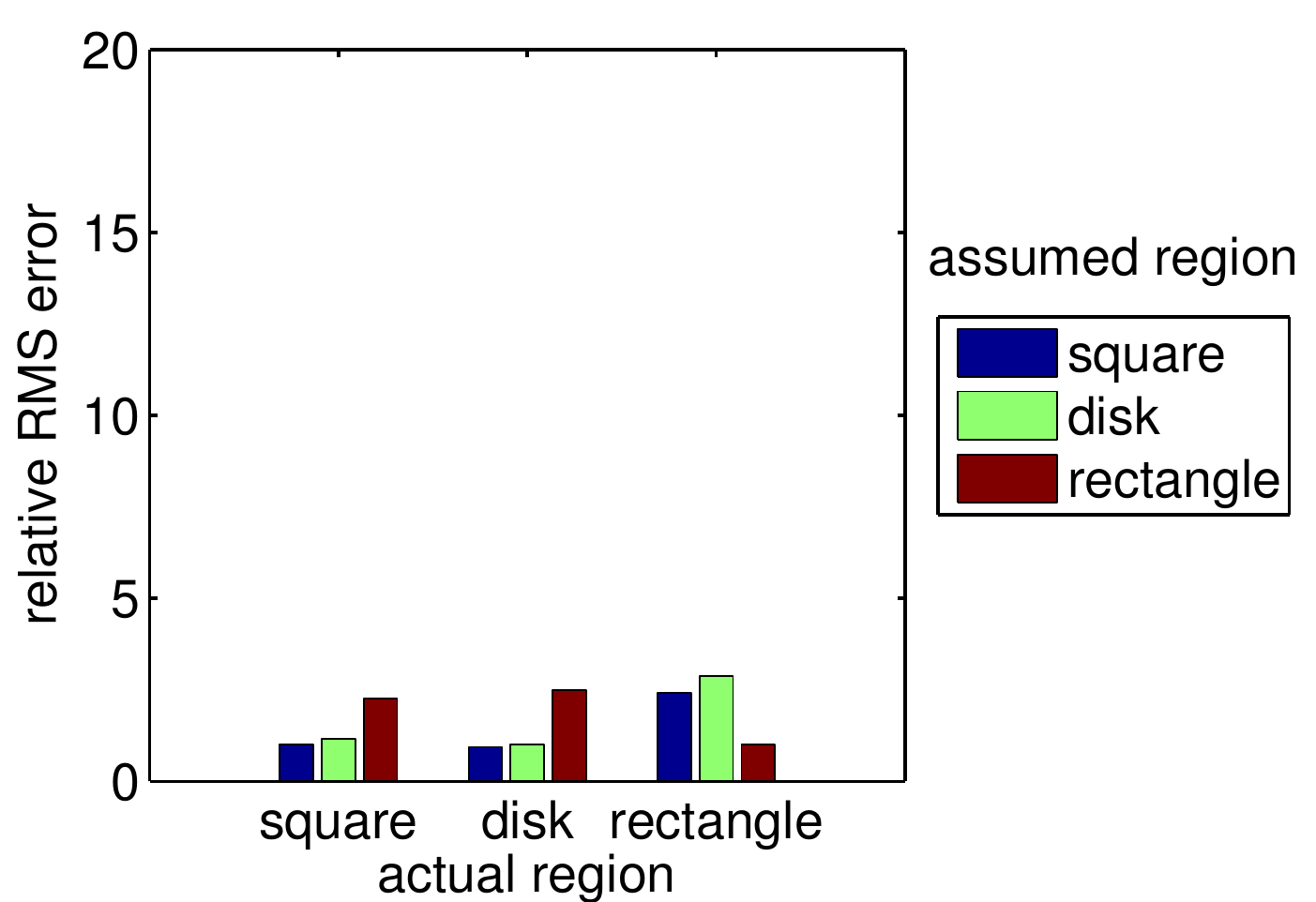}
      \caption{(b) $s = 20.0$.}
    \end{subfigure}


    \caption{Relative RMS error when the assumed region shape is
      incorrect. Note that errors are much smaller for large $s$
      values, because boundary affects play a lesser role. }
    \label{fig:wrong_shape}
  \end{center}
\end{figure}

More importantly, it is possible that we will not know the exact shape
of the region of interest, or that the shape used in estimation is an
approximation of an irregular region. In these cases we could use the
empirical estimates of the potential links, but if the shape of the
region is unknown, it is important to consider the impact of an
erroneous decision. \autoref{fig:wrong_shape} shows that impact by
highlighting the relative size of the induced error as a result of
assuming the wrong shaped region for three different region shapes (a
square, a disk and a rectangle with an aspect ratio of 2:1).

Figure \ref{fig:wrong_shape}(a) shows that there is a certain
robustness, as choosing the disk instead of the square or visa versa
leads to only a small increase in error. On the other hand, assuming a
square or disk when the region is really a rectangle causes much
larger errors. The difference is that the longer thinner rectangle has
a considerably larger boundary to area proportion, and the boundary effects
play a larger role in the random line distribution.

Comparison of \autoref{fig:wrong_shape}(a) and (b) suggests that the
effect is more significant for smaller values of $s$. This is because
larger values of $s$ result in shorter links on average hence boundary
effects play a smaller role. So unsurprisingly, the estimates are more
robust for large $s$.

\begin{figure}[tbp]
  \begin{center}
    \includegraphics[width=\figurewidthA]{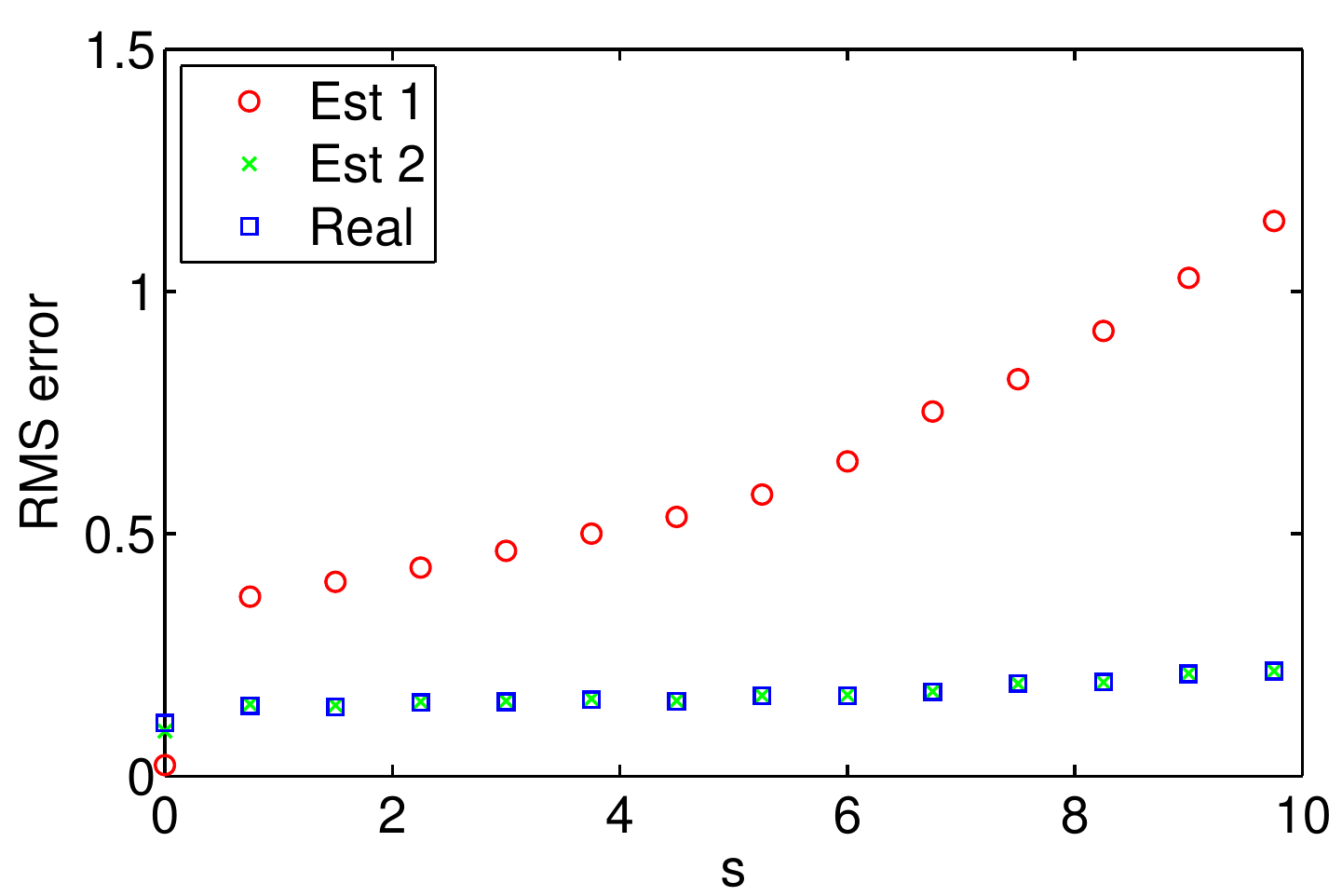}
    \caption{The impact of estimating the size of the region. Real is
      the MLE given the real region size. Est 1 is the MLE given the
      region size is estimated using the longest link, and Est 2 is
      the MLE given the region size is estimated from the longest
      inter-node distance. We see that estimating the region size
      using the inter-node distances introduced minimal errors, but
      using the link distances increases the errors, except for very
      small $s$. \label{fig:est_L}}
  \end{center}
\end{figure}

We also consider the impact of estimating the size of the
region. \autoref{fig:est_L} shows RMS errors for a disk-shaped region
when the diameter is known (Real), when it is estimated by taking
the largest link-distance (Est 1), and when it is estimated using the
largest inter-node distance (Est 2). 

The loss in estimator accuracy when the size of the region is unknown
is clear. Further, the inaccuracy increases with $s$ because as
$s$ increase, the longest observed link is less likely to be a proxy
for the actual size of the region. We could correct for this bias, but
instead we compare to estimates of the region size formed from the
maximum inter-node distance, and we see that there is
almost no loss of accuracy in the estimator. 

The conclusion is that some knowledge of the region over which the graph is
drawn is important. Very poor estimates of the region can result in
much worse estimator accuracy, but there is a fair amount of
robustness so the estimate of the region shape and scale 
need not be perfect and estimates drawn from the data are sufficient. 

\begin{figure}[tbp]
  \begin{center}
    \includegraphics[width=\figurewidthA]{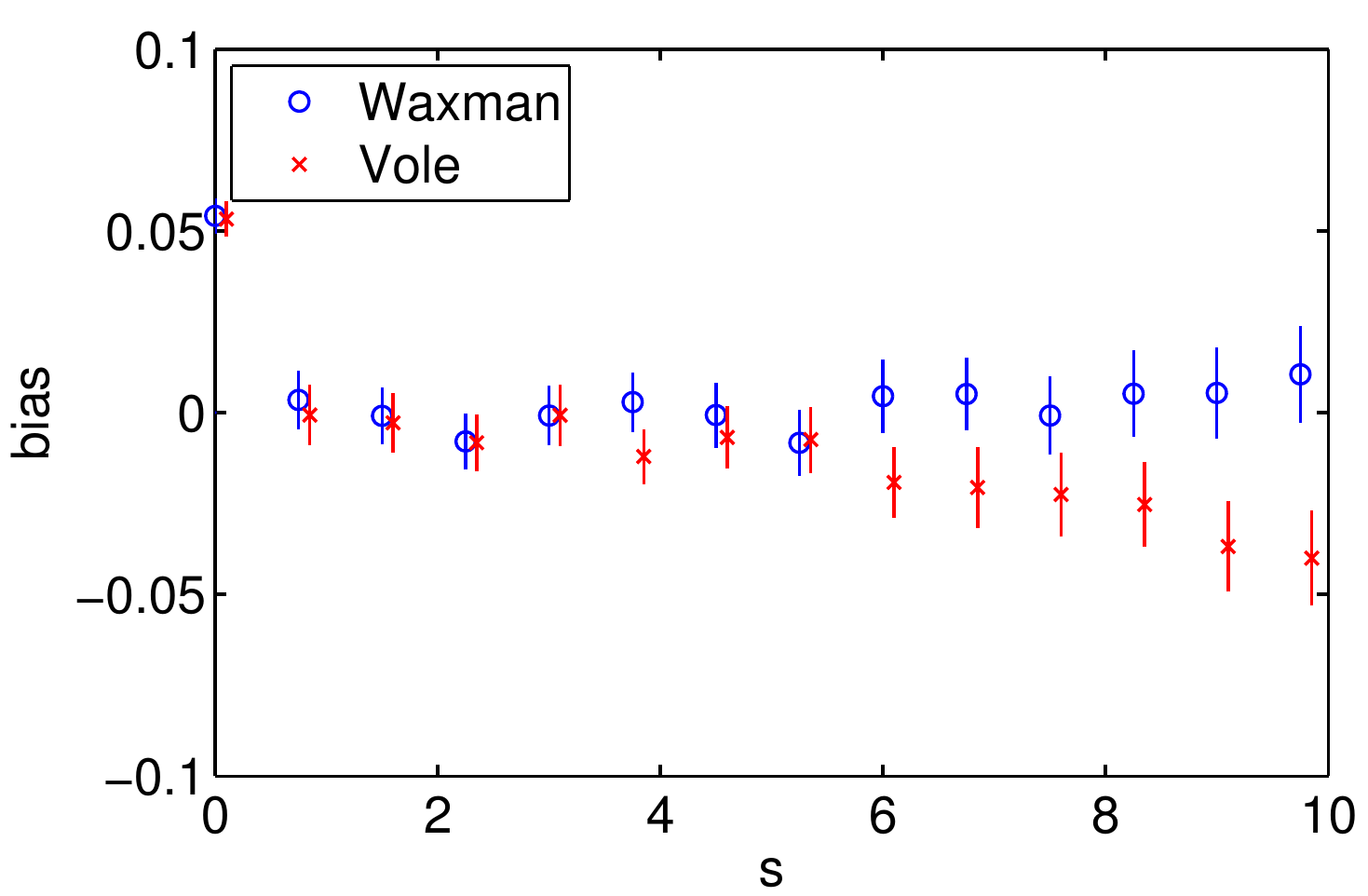}
    \caption{Estimates $\hat{s}$ when the underlying network is not
      Waxman, but rather drawn from the DASTB vole-contact model. For
      smaller $s$, using the wrong model introduces minimal addition
      error.}
    \label{fig:est_vole}
  \end{center}
\end{figure}

Finally, to test how robust the model outputs are when the underlying model is
incorrect, \ie the distance deterrence function is not Waxman, 
we test the performance of the MLE when the underlying graph is
actually the DASTB model~\citep{davis14:_spatial} used for the vole
contact graph data. We use this model as a comparison, because it is
very similar in concept, but has different details.
\autoref{fig:est_vole} shows that there is some very small induced
bias in the estimates as a result of applying the model in the wrong
situation, but these would have almost no influence on the
conclusions draw from $\hat{s}$ estimates. 

\subsection{Estimating $\hat{q}$}

In most of the preceding work we have only considered the accuracy of
the estimates for $s$. The estimated accuracy of $q$ is that of the
binomial model parameter estimate which is well understood   
except for the manner in which errors in $\hat{s}$ propagate into the
estimate $\hat{q}$ through $\tilde{G}(\hat{s})$. 

\begin{figure}[tbp]
  \begin{center}
    \includegraphics[width=\figurewidthA]{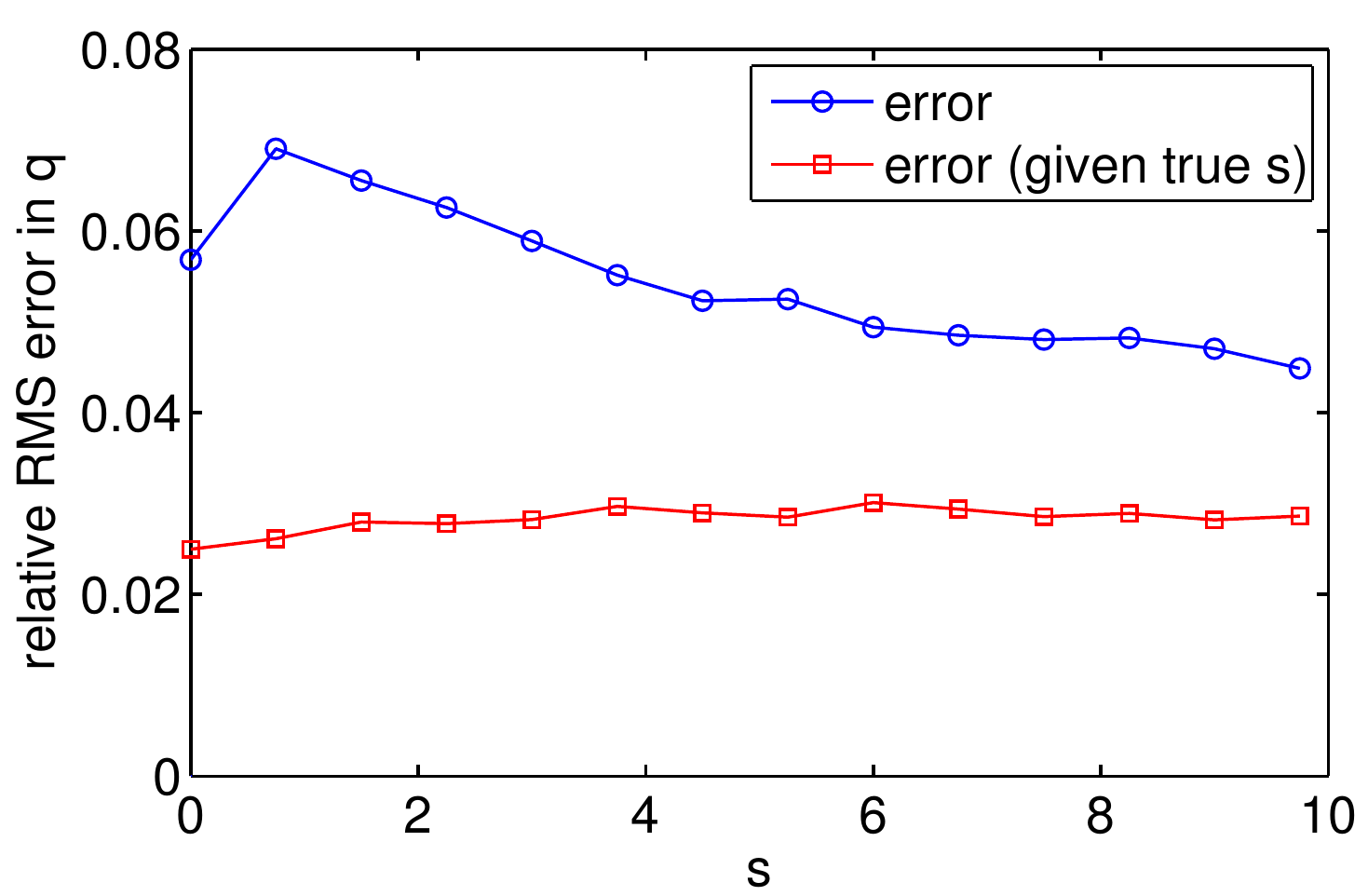}
    \caption{Relative error for the MLE $\hat{q}$ ($n=1000$,
      $\bar{k}=3$). Squares indicate the RMS error relative to $q$
      assuming we were given the true $s$ values, and circles indicate
      the error of the actual estimator.}
    \label{fig:q}
  \end{center}
\end{figure}

\autoref{fig:q} shows the RMS errors in $\hat{q}$ relative to the size
of $q$, given the estimated and true values of $s$.  We can see from
the figure that the errors in $\hat{q}$ are roughly doubled by the
uncertainty in $\hat{s}$, so this error does have a significant
effect. However, the errors in $\hat{q}$ are still small: \ie of the
order of 4\%.

\section{Case studies}
 
The previous section looked at accuracy on simulated data, where we
know the ground truth. This section considers how well the method
works on real data. We consider three datasets, one biological and two Internet.

\subsection{Vole data}

\begin{figure}[tbp]
  \begin{center} 
    \begin{subfigure}[t]{\figurewidthB}
      \includegraphics[width=\textwidth]{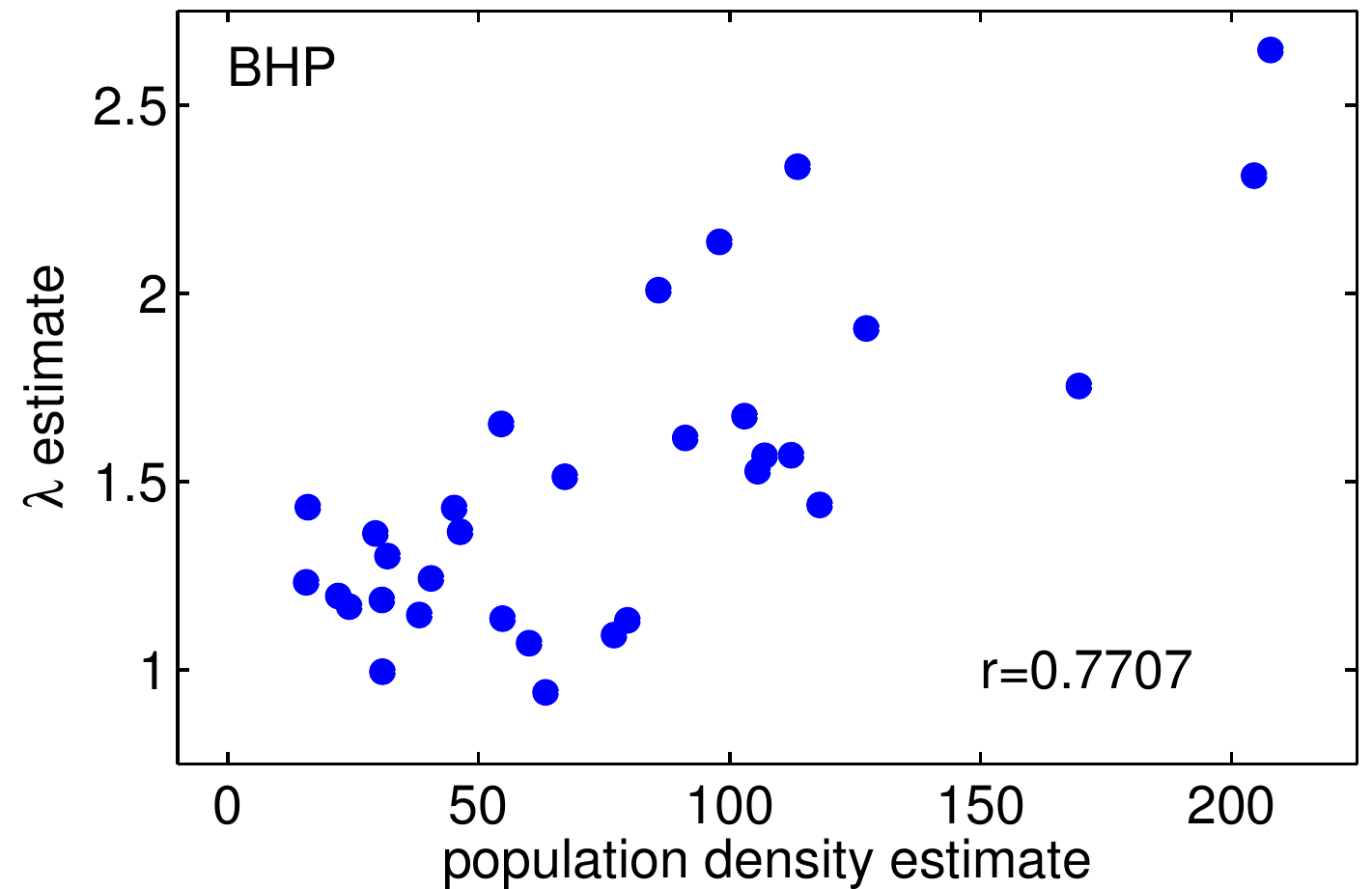}
      \caption{(a) Reproduction of Figure 5(a) from
        \citep{davis14:_spatial}.}
    \end{subfigure}
    \hfil
    \begin{subfigure}[t]{\figurewidthB}
      \includegraphics[width=\textwidth]{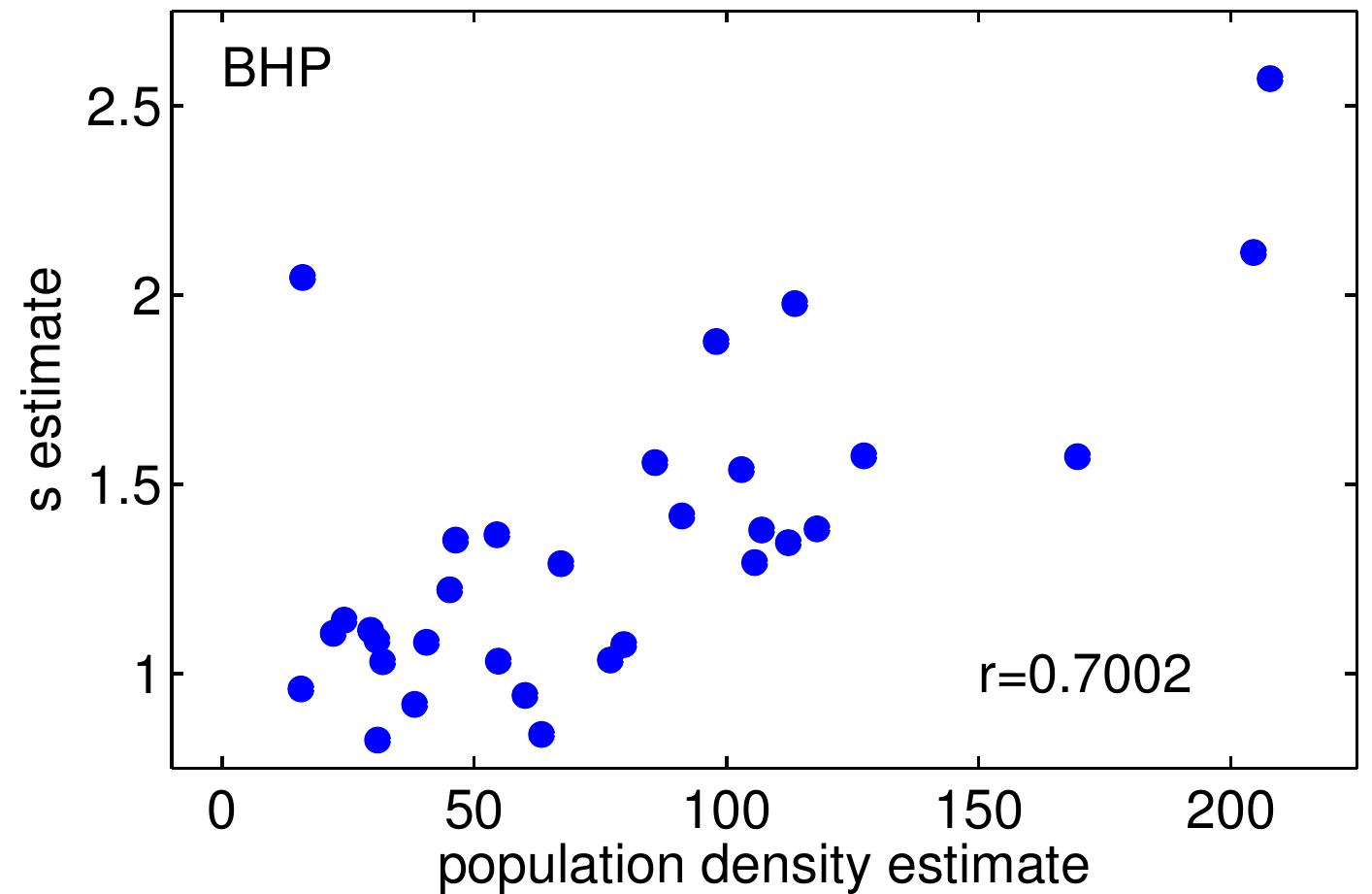}
      \caption{(b) Equivalent plot of $\hat{s}$ estimates against
        population density.}
    \end{subfigure}
    \caption{Vole data: data drawn from the Black Blake Hope (BHP)
      site: $r$ reports the correlation coefficients. \
      \label{fig:vole_data}}
  \end{center}
\end{figure}

\citet{davis14:_spatial} collected a large mark-recapture dataset on
the English–Scottish border over a 7-year period to study
relationships between field voles {\em (Microtus agrestis)}. In
particular, they considered a number of models for the contact graphs
created by relationships between the voles.
A grid of trap locations was created on each of four sites and over
64 capture periods traps were emptied multiple times. Contact graphs
were built by associating voles that were caught in the same trap on
separate occasions. In their main dataset, the time periods were
divided into pairs and each pair used to construct one graph, so in
total 32 graphs were created for each site. One might argue with the
methodological accuracy of constructing contact graphs from trap data
in this manner, but there is no doubt that this is an unusually
large dataset of related graphs. It is rare to have more than a few
example contact graphs in any one application. The large number of
these graphs allows some views of estimates that are not usually
available.

\citet{davis14:_spatial} fitted multiple models to the data including
a GER random graph, and their DASTB model, which resembles the Waxman
graph. Their model incorporated the capture process as an explicit
component and they form an edge in the graph when one or more contacts
occurred. 

Using the original raw data provided by the authors we reconstructed
their contact graphs and confirmed we could reproduce the results
obtained in \citep{davis14:_spatial}.  Estimates of $\lambda$ (their
distance parameter) were recomputed using the complimentary log-log
regression to confirm that the recreated graphs and their measured
statistics were correct, and that the units were
correct. \autoref{fig:vole_data}(a) reproduces \citet[Figure
5(a)]{davis14:_spatial}, demonstrating that the data has been
correctly reconstructed.

Next we performed our own estimates using the MLE with the underlying
region being a square of side 10 (comprised of 10 traps with a 5m
spread between each, thus matching the
experiment. \autoref{fig:vole_data}(b) shows our estimates (we have
repeated the same procedure for all four original datasets with
similar results). We can see from the graph that although estimates of
$s$ are not the same as $\lambda$, the general conclusion of
\citep{davis14:_spatial} that there is a correlation between
population density and the spatial contact parameter is supported
with almost the same correlation coefficient being reported for both
datasets.

The results emphasise the original findings of Davis~\etal, notably that:
\begin{itemize}

\item Animal contact graphs are well modelled by spatial graphs (in
  our case there is a small difference in the particular model, but
  the general behaviour is similar).

\item This modelling provides a means of quantifying differences between two
  groups of animals (separated in time or space).

\item The parameter $\hat{s}$ tells us something important that other
  many models omit, namely that the distance voles travel when
  contacting each other decreases with population density. Thus has a
  profound importance for disease transmission \citep{davis14:_spatial}
  in that population mixing may not increase with population density
  as quickly as is often assumed.

\end{itemize}
More generally, models can be used for the usual suite of
purposes such as simulation and extrapolation. For example,
once we know the length scale at which contact is
unlikely, we might use this in the design of capture/recapture
experiments by choosing the scale and resolution of the capture setup.
  
\subsection{Internet dataset 1} 
 
\citet{Lakhina:2002:GLI:637201.637240} undertook one of the first
attempts to formally quantify the exponential decrease of link
likelihood as a function of distance. The authors compared two sets of
data and found consistent results between them. They provided one of
these datasets (the Mercator data) to us for comparison.
Lakhina~\etal\ separated the data into three regions, and analysed
these separately. \autoref{tab:mercator_summary} provides a brief
summary of the three regions.

We do not argue that this network is random in any real sense:
in fact the Internet networks are the result of design.  However,
fitting a Waxman-like graph to these is instructive in that it shows
how engineering constraints lead to distance-sensitive link placement.
 
We can see first from the data that the graphs are large but very
sparse, with average node degrees of around 2. In this dataset we
also have node locations so we can derive distances between all pairs
of nodes and thus apply all of the possible techniques considered here.
However, it is not feasible to use the GLM approach for this scale of problem.

\citet{Lakhina:2002:GLI:637201.637240} applied log-linear regression
to the question. We applied the MLE and compared the results to
those found in their work. \autoref{tab:mercator_results} provides a
comparison between various estimates, including the original values
reported in \citet{Lakhina:2002:GLI:637201.637240} given in the second
column under $\hat{s}_{Lakhina}$. Units are {\em per 1000 miles} ---
we use Imperial units to be consistent with the original paper. The
third column provides our equivalent estimate. There are small
differences, presumably because of differences in the exact
numerical procedures applied.

\begin{table}
  \caption{Summary of Mercator datasets: number of nodes, edges, and
    the average link distance.\label{tab:mercator_summary}}
  \centering
  \def\arraystretch{1.2}    
%
\begin{tabular}{rrrrr}
  \hline
  Region & n & e & $\bar{d}$ (miles)  \\ 
  \hline
     USA & 123426 & 152602 &    384.7  \\ 
  EUROPE &  32928 &  30049 &    319.5  \\ 
   JAPAN &  14318 &  16665 &    317.6  \\ 
  \hline
\end{tabular}

\end{table}

The fourth column of the table shows the MLE-E values estimated for the
datasets. We use the Empirical estimator because the region
shapes are irregular (\eg the USA), and we want to avoid approximation
errors arising from the region shape. 

We see considerable discrepancies which are larger than can be
explained by errors in the log-linear regression approach. However,
reading \citet{Lakhina:2002:GLI:637201.637240}, we can see that their
estimates are over a truncated range of distances for two reasons: 
\begin{itemize}
\item The node locations they use are artificially quantised by the
  Geolocation procedure so some nodes appear to have exactly
  the same position, and hence zero distance, when actually there is a
  positive distance between the nodes; and 

\item They found that the exponential distance-deterrence function fit
  the data only up to some threshold distance.

\end{itemize}
In their (and our comparison) log-linear regressions the range over
which we perform the regression is restricted to be between these
bounds.

In order to provide a fair comparison we also modified the MLE-E by
censoring the potential edges used in forming the CDF and in computing
the average edges distance.  \autoref{tab:mercator_results} shows the
results under $s_{MLE-T}$ to be closer to being consistent with the
log-linear regression.
  
The results point to one valuable feature of the log-linear
regression, which is that it comes with diagnostics. Examination of
the fit indicates whether the model is appropriate or not. The MLE
requires additional effort to provide similar diagnostics. On the
other hand, there are significant issues with the log-linear
regression. Apart from being less accurate, there is the question of
bin size which simple experiments seem to suggest has more effect on
estimates than one might hope.

Ultimately, all of the methods suggest strongly that a spatial
component should be part of any model for Internet linkages. This is
entirely consistent with the intuition of engineers who work on such
networks: long links cost more, and so are  rarer.
 
\subsection{Internet dataset 2} 

\begin{figure}[tbp]
  \begin{center}
    \includegraphics[width=\figurewidthA]{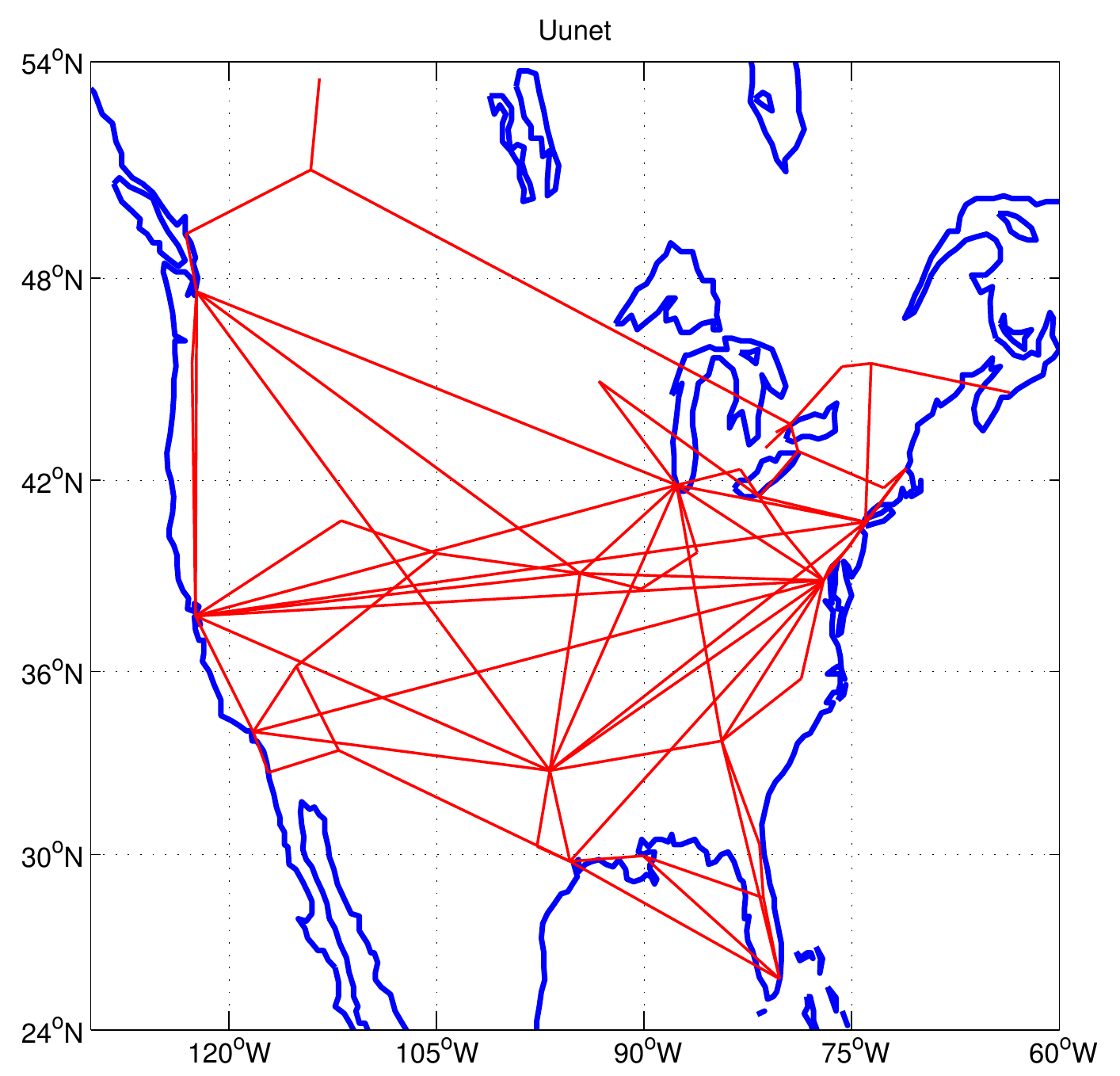}
    \caption{The Uunet network from the Internet Topology Zoo.}
    \label{fig:uunet}
  \end{center}
\end{figure}

Finally we apply the MLE to a set of real networks from the Internet
Topology Zoo~\citep{Zoo}. The networks -- taken at Point of Presence
(PoP) level -- show the connectivity of a sample of seven major
network operators in different regions of the world.
\autoref{fig:uunet} shows one example, and the seven are summarised and results shown in
\autoref{tab:zoo_results}. 

\begin{table}
  \caption{Estimates: $\hat{s}_{Lakhina}$ are the values from
    \citet{Lakhina:2002:GLI:637201.637240}; and $\hat{s}_{log-linear}$
    are our corresponding estimates. Units are {\em per 1000 miles}.
    The $\hat{s}_{MLE-E}$ values are derived from our 
    Empirical MLE, and the$\hat{s}_{MLE-T}$ values from a version of
    the MLE  with distance data truncated in the same manner as the original log-linear
    estimates. 
    \label{tab:mercator_results}}
  \centering
  \def\arraystretch{1.2}  
%
\begin{tabular}{rrrrr}
  \hline
  Region & $\hat{s}_{Lakhina} $& $\hat{s}_{log-linear}$ & $\hat{s}_{MLE-E}$ & $\hat{s}_{MLE-T}$  \\ 
  \hline
      USA &     6.91 &     6.38 &     2.75 &     6.63 \\ 
      EUROPE &    12.80 &    12.81 &    30.92 &    10.09 \\ 
      JAPAN &     6.89 &     6.71 &    45.91 &     7.30 \\ 
  \hline
\end{tabular}

\end{table}  

These graphs are
much smaller as a result of being drawn at the PoP level, and so
errors in the estimates will be correspondingly larger. However, there
is a very clear variation in the parameters. 

Firstly, although the units for $s$ here are per 1000 km, we can
clearly see that the $s$ estimates are a good deal smaller than for
the previous dataset, largely because of looking at a
single network at a time, at the PoP level. Router-level networks
(such as the previous dataset) include numerous routers, many of them
connected, in the same cities. The interconnects between networks also
tend to be shorter as they are mediated through Internet exchange
points or carrier hotels.  However, at the PoP level, all of these
small-scale details are invisible, and consequently we see less strong
distance deterrence.

Secondly, we see a wide variety of $s$ values, with little
relationship to region. Networks are just different. They are built
with different goals, and different cost structures and
constraints. Some are new, and others are older and incorporate legacy
components.

The wide variety of network structures has been observed in the Zoo
data before \citep{bowden14:_cold}, but not with respect to spatial
structure. This reinforces the message that estimates of model
parameters such as $s$ provide alternative methods to compare network
behaviour. 

\begin{table}
  \caption{Estimates for Internet Topology Zoo data. Units for $s$ are
  {\em per 1000 km}. \label{tab:zoo_results}}
  \centering
  \def\arraystretch{1.2}  
%
\begin{tabular}{rr|rrr|rr}
  \hline
  Network & Region & $n$ & $e$ & $\bar{d}$ & $\hat{s}_{MLE-E}$ \\ 
  \hline
        Aarnet &    Australia & 19 & 24 & 695.6 & 2.20  \\
         Iinet &    Australia & 9 & 12 & 1472.9 & 0.30  \\
      BtEurope &       Europe & 22 & 35 & 606.3 & 2.05  \\
          Colt &       Europe & 153 & 177 & 160.2 & 9.48  \\
       Abilene &          USA & 11 & 14 & 1007.0 & 1.39  \\
   Internetmci &          USA & 19 & 33 & 927.7 & 1.25  \\
         Uunet &          USA & 42 & 77 & 966.9 & 1.09  \\
  \hline
\end{tabular}

\end{table} 

\section{Discussion and Conclusion}

This paper presents the MLE for the parameters of the Waxman graph and
demonstrates its accuracy in comparison to alternative estimators.
The MLE has two advantages. Firstly it can guarantee $O(n)$
computational time complexity and constant memory usage by using only
a sample of the edges that exist in a graph to estimate $s$. Secondly
it can be applied in domains where the coordinates of nodes are
unknown and/or edge lengths may be weights in some arbitrary process.
Using the MLE we have shown that real networks have a considerably
wider range of $s$ parameters than is typically used in the
literature.

The next question might be: ``Is the Waxman model a better model than
X?''  where X might be the GER random graph, or some other model.
\citet{davis14:_spatial} considered this question using AIC, but for
the simple question of testing whether there should be a distance term
or not, it is perhaps more natural to use a hypothesis test. We applied
the standard likelihood ratio test and found that standard threshold
choices resulted in poor ability to specify the significance. It
seems the problem arises from the failure of the independence
assumption. If so the problem is correctable, but we leave the
question of correctly determining the cutoff threshold for future
work.

\section{Acknowledgements}

We would like to thank \citet{davis14:_spatial} for providing us with
their vole data, and \citet{Lakhina:2002:GLI:637201.637240} for
providing the first Internet dataset. The Zoo data is publicly
available at \url{www.topology-zoo.org}.

This work was supported by the Australian Research Council through
grant DP110103505, and the Centre of Excellence for Mathematical \&
Statistical Frontiers.

{\footnotesize 
\setlength{\parskip}{-1mm}
\bibliographystyle{rss}
\providecommand{\urlprefix}{Online: }
\bibliography{paper}
}

\end{document}